\documentclass[10pt]{amsart}

\usepackage[T1]{fontenc}
\usepackage{lmodern}       

\usepackage{microtype}

\usepackage{amsmath}
\usepackage{amssymb}
\usepackage{amsthm}
\usepackage{amsxtra}
\usepackage{mathtools}
\usepackage{mathrsfs}
\usepackage{stmaryrd}
\usepackage{leftindex}
\usepackage{tensor}
\usepackage{mhchem}
\usepackage{colonequals}
\numberwithin{equation}{section}
\usepackage[draft]{graphicx}    
\usepackage{wrapfig}
\usepackage{float}
\usepackage{caption}
\usepackage{subcaption}
\usepackage{epsfig}
\usepackage{tikz-cd}

\usepackage{array}
\usepackage{booktabs}
\usepackage{multirow}
\usepackage{multicol}
\usepackage{diagbox}
\usepackage{numprint}
\usepackage{enumerate}
\usepackage{enumitem}
\setlist{itemsep=4pt}
\usepackage{fullpage}

\usepackage{soul}
\usepackage{comment}
\usepackage{xcolor}

\definecolor{webcolor}{rgb}{0.8,0,0.2}
\definecolor{webbrown}{rgb}{.6,0,0}

\usepackage{hyperref}
\hypersetup{colorlinks=true,urlcolor=blue,citecolor=blue,linkcolor=blue}
\usepackage[alphabetic,backrefs,lite]{amsrefs}

\theoremstyle{definition}
\newtheorem{thm}{Theorem}[section]
\newtheorem{lemma}[thm]{Lemma}
\newtheorem{cor}[thm]{Corollary}
\newtheorem{prop}[thm]{Proposition}

\newtheorem{program}[thm]{Program}
\newtheorem*{programb}{Mazur's Program B}

\theoremstyle{definition}
\newtheorem{defn}[thm]{Definition}

\newtheorem{algorithm}[thm]{Algorithm}

\theoremstyle{remark}
\newtheorem{remark}[thm]{Remark}

\newcommand{\defi}[1]{\textsf{#1}}     

\newenvironment{romanenum}
  {\hfill \begin{enumerate}[label=({\roman*})]}
  {\end{enumerate}}

\renewcommand{\AA}{\mathbb A}

\newcommand{\CC}{\mathbb C}

\newcommand{\NN}{\mathbb N}
\newcommand{\PP}{\mathbb P}
\newcommand{\QQ}{\mathbb Q}
\newcommand{\RR}{\mathbb R}

\newcommand{\ZZ}{\mathbb Z}

\newcommand{\Zhat}{\widehat\ZZ}
\newcommand{\Ztimes}{\widehat\ZZ^{\times}}

\newcommand{\calD}{\mathcal D}

\newcommand{\calF}{\mathcal F}
\newcommand{\calG}{\mathcal G}

\newcommand{\HH}{\mathcal H}

\newcommand{\calL}{\mathcal L}

\newcommand{\calR}{\mathcal R}

\newcommand{\calX}{\mathcal X}

\newcommand{\scrA}{\mathscr A}

\newcommand{\scrF}{\mathscr F}

\newcommand{\scrL}{\mathscr L}

\newcommand{\GL}{\operatorname{GL}}
\newcommand{\SL}{\operatorname{SL}}

\newcommand{\PSL}{\operatorname{PSL}}

\newcommand{\coef}{\operatorname{coef}}

\newcommand{\Gal}{\operatorname{Gal}}

\newcommand{\ord}{\operatorname{ord}}

\def\bbar#1{\setbox0=\hbox{$#1$}\dimen0=.2\ht0
  \kern\dimen0 \overline{\kern-\dimen0 #1}}
\newcommand{\Qbar}{\overline{\mathbb Q}}

\newcommand{\ddfrac}[2]{\frac{\displaystyle #1}{\displaystyle #2}}

\newcommand{\Ssl}[1]{\text{SL}_{2}(#1)}
\newcommand{\Ggl}[1]{\text{GL}_{2}(#1)}

\newcommand{\commentinv}[1]{{}}

\usepackage[euler-digits,small]{eulervm}
\makeatletter
\renewcommand{\[}{\begin{equation*}}
\renewcommand{\]}{\end{equation*}}
\makeatother
\begin{document}
\title[A Classification of Low Genus Modular Curves]{A Classification of Low Genus Modular Curves}
\author{Eray Karabiyik}
\address{Department of Mathematics, Cornell University, Ithaca, NY 14853, USA}
\email{ek693@cornell.edu}



\begin{abstract} 

Let $G$ be an open subgroup of $\GL_2(\Zhat)$ satisfying $\operatorname{det}(G)=\Zhat^{\times}$ and $-I\in G$. Associated to $G$, there is a modular curve $X_G$ defined over $\QQ$ which parametrizes elliptic curves with $G$-level structure. We give a classification of modular curves with genus equal to a fixed nonnegative integer $g$. In particular, we show that all modular curves of genus $g$ lie in finitely many families of $\mathbb{Q}^{\operatorname{ab}}$-twists of modular curves. We also describe an algorithm for computing all families of modular curves of a fixed genus $g$ and use this to compute projective models for these modular curves. This algorithm has been fully implemented for $g \leq 12$.

In an appendix with Yongyuan Huang and Rakvi, we describe two algorithms, one for computing rational functions on modular curves given as quotients of two modular forms of the same weight, and the other for computing models for the generic fibers of universal elliptic curves $\mathcal{E}_G$, where $G\subset \GL_2(\Zhat)$ is an open subgroup satisfying $\operatorname{det}(G)=\Zhat^{\times}$ and $-I\notin G$.

\end{abstract}

\maketitle
\section{Introduction}  \label{S:introduction}

\indent Let $E$ be a non-CM elliptic curve defined over the rational numbers. Fix an algebraic closure $\overline{\QQ}$ of $\QQ$. Let $\Gal_\QQ:=\Gal(\Qbar/\QQ)$ denote the absolute Galois group of $\QQ$. For any positive integer $N$, let $E[N]$ be the $N$-torsion of $E(\overline{\QQ})$, it is a free $(\ZZ/N\ZZ)-$module of rank $2$. The group $\Gal_\QQ$ acts naturally on $E[N]$ and respects the group structure. This gives rise to a Galois representation
\[
\rho_{E,N}\colon \Gal_\QQ \rightarrow \operatorname{Aut}(E[N]) \cong \operatorname{GL}_2(\ZZ/N\ZZ).
\]
\indent Fixing compatible bases for $E[N]$ for all $N\geq 1$, and taking the inverse limit, one gets the adelic representation
\[
\rho_{E}\colon \Gal_\QQ \longrightarrow \operatorname{Aut(E_{\operatorname{tors}})} \cong \operatorname{GL}_2(\Zhat)
\]
 where $\Zhat$ is the profinite completion of $\ZZ$. The image of $\rho_E$ is uniquely determined up to conjugacy in $\GL_2(\Zhat)$. The group $\rho_E(\Gal_\QQ)$ is a closed subgroup of $\operatorname{GL}_2(\Zhat)$ with respect to the profinite topology. In \cite{Serre-Inv72}, Serre proved that $\rho_E(\Gal_\QQ)$ is an open subgroup of $\GL_2(\Zhat)$.
 
\indent Let $\chi_{\operatorname{cyc}}:\Gal_\QQ \rightarrow \Zhat^{\times}$ be the cyclotomic character. Using the Weil pairing on $E$, one can show that $\operatorname{det}\circ \rho_E$ agrees with $\chi_{\operatorname{cyc}}$, and hence the image of $\rho_E$ has full determinant, i.e, $\det(\rho_E(\Gal_\QQ))=\Zhat^{\times}$. Let $G$ be an open subgroup of $\operatorname{GL}_2(\Zhat)$ such that $\operatorname{det}(G)=\Zhat^{\times}$ and $-I\in G$. Associated to $G$, there is a smooth, projective, geometrically integral curve $X_G$ over $\QQ$ that parametrizes elliptic curves with $G$-structure which will be explicitly defined in \S \ref{S:Background}. A \defi{modular curve} is any curve of the form $X_G$.

\indent In this paper, we classify modular curves of genus $g$ into finitely many families of $\mathbb{Q}^{\operatorname{ab}}$-twists. We state here a simplified version; it is proven as part of Theorem \ref{T:main families1}.

\begin{thm}[Classification of modular curves of fixed genus]\label{T:ilk}
    Fix a non-negative integer $g$. There are finitely many sets (called \defi{families of modular curves}) $(\mathscr{F}_i)_{i\in I}$ of modular curves such that;

\begin{enumerate}
    \item The family $\mathscr{F}_i$ consists of twists of a single modular curve $X_{G_i}$ for all $i\in I$,
    \item If $X_G$ is a modular curve of genus $g$, then $X_G\in \mathscr{F}_i$ for some $i\in I$.
\end{enumerate}
\end{thm}

Using these families, we also describe an algorithm that efficiently computes projective models of modular curves of genus $g$. We have computed all such families of modular curves of genus at most $24$. For genus at most $12$, we have precomputed representative models and implemented the algorithm. A \texttt{Magma} \cite{Magma} package implementing the algorithm can be found at \cite{repo}.  

Our \texttt{Magma} package has been used by \cite{LMFDB} to compute projective models for more than one million modular curves in their modular curves database.




\subsection{Modular curves}
Let $G\subseteq \operatorname{GL}_2(\Zhat)$ be an open subgroup such that $\operatorname{det}(G)=\Zhat^{\times}$ and $-I\in G$. We will define the associated modular curve $X_G$ in \S \ref{S:Background}. It is a smooth, projective, geometrically irreducible curve defined over $\QQ$. An inclusion $G\subseteq G' \subseteq \GL_2(\Zhat)$ induces a nonconstant morphism of curves $ X_G\rightarrow X_{G'}$.
For the group $\operatorname{GL}_2(\Zhat)$, we have particular isomorphism $X_{\operatorname{GL}_2(\Zhat)}\cong\PP_{\QQ}^1=\AA_\QQ^1 \cup \{\infty\}$, and call $X_{\operatorname{GL}_2(\Zhat)}$ the \defi{j-line}. Taking $G'=\GL_2(\Zhat)$, we have the associated $j$-map 
\[
\pi \colon X_G \longrightarrow \PP_\QQ^1.
\] 
\indent Let $\rho_E^{\ast}\colon \Gal_\QQ \rightarrow \GL_2(\Zhat)$ be the dual representation of $\rho_E$ defined by $\rho_E^{\ast}(\sigma)=\rho_E(\sigma^{-1})^\intercal$. The curve $X_G$ has the following property, cf. \cite{explicitopenimage}*{Proposition 3.5}. 
\begin{prop}  \label{P:intial moduli property}
Let $E$ be any elliptic curve defined over $\QQ$ with j-invariant $j_E\notin \{0,1728\}$.  Then $\rho_{E}^*(\Gal_\QQ)$ is conjugate in $\GL_2(\Zhat)$ to a subgroup of $G$ if and only if $j_E$ is an element of $\pi_G(X_G(\QQ)) \subseteq \QQ \cup \{\infty\}$.  
\end{prop}
Throughout the paper, by the \defi{genus} of an open subgroup $G\subseteq \GL_2(\Zhat)$, we mean the genus of the associated modular curve $X_G$. If $G_1$ and $G_2$ are open subgroups of $\GL_2(\Zhat)$ that are conjugate to each other, then associated modular curves are isomorphic; so we only consider open subgroups up to conjugacy in $\GL_2(\Zhat)$.  The \defi{level} of an open subgroup $G$ of $\operatorname{GL}_2(\Zhat)$ is the smallest positive integer $N$ such that $G$ contains the kernel of the reduction modulo $N$ map $\operatorname{GL}_2(\Zhat) \rightarrow \operatorname{GL}_2(\ZZ/N\ZZ)$.

For classical families of modular curves ($X_0(N)$, $X_1(N)$, $X_{ns}(N)$, $X_{ns}^+(N)$, $X_{s}(N)$ and so on) there are many algorithms to compute models in the literature. On the other hand, for an arbitrary open subgroup $G\subseteq\GL_2(\Zhat)$ with $-I\in G$ and $\det(G)=\Zhat^{\times}$, there are limited methods to compute a model for $X_G$. One such algorithm is the one implemented in \cite{explicitopenimage}, see \S \ref{S:Background}.  Let $g\geq0$ be an integer. In \S \ref{S:algorithm}, using a twisting argument, we describe a related algorithm to compute the model of any modular curve $X_G$, whose genus is $g$. When $X_G$ is a geometrically non-hyperelliptic modular curve of genus at least $2$, the algorithm computes the image of the canonical map as the projective model.\\
\indent There are many equivalent definitions of the modular curve $X_G$. One can define $X_G$ by explicitly defining its function field or as the generic fiber of the coarse moduli space of the stack $M_G$ over $\ZZ[1/N]$ parametrizing generalized elliptic curves with $G$-level structure, see \cite{MR0337993} for details. One can also refer to \cite{katzmazur} for the fine arithmetic of modular curves, where the level structure has a meaning over schemes where $N$ is not invertible.   We will opt to define $X_G$ through a certain graded ring of modular forms which is defined in \S \ref{S:Background}.

\subsection{Agreeable groups}
\indent Fix a non-negative integer $g$. Let $G$ be an open subgroup of $\operatorname{GL}_2(\Zhat)$ that has full determinant and contains $-I$. For our classification of modular curves of genus $g$, we introduce a special kind of subgroup of $\operatorname{GL}_2(\Zhat)$ that contains the group $G$.

\begin{defn}
    An open subgroup $H \subseteq \operatorname{GL}_2(\Zhat)$ is \defi{agreeable} if $\operatorname{det}(H)=\Zhat^{\times}$, $H$ contains the scalar matrices, i.e. $\Zhat^{\times}\cdot I \subseteq H$, and the levels of $H$ in $\operatorname{GL}_2(\Zhat)$ and $H \cap \operatorname{SL}_2(\Zhat)$ in $\operatorname{SL}_2(\Zhat)$ have the same odd prime divisors.
\end{defn}
\indent There exists a group $G\subseteq \calG \subseteq \GL_2(\Zhat)$, called the \defi{agreeable closure} of $G$, which is minimal, with respect to inclusion, among all agreeable subgroups that contain $G$. The group $G$ is normal in $\calG$ and satisfies $[G,G]=[\calG,\calG]$ i.e., their commutator subgroups agree. Since   there is a nonconstant morphism  $X_G\rightarrow X_\calG$, the genus of $\calG$ is less than or equal to the genus of $G$.

The set of agreeable subgroups is closed under conjugation in $\operatorname{GL}_2(\Zhat)$. We show in \S\ref{S:finiteness} that there are only finitely many agreeable subgroups of genus at most $g$.

\subsection{Families attached to a pair}\label{SS:famsfirst}
\indent Let $\calG \subseteq \operatorname{GL}_2(\Zhat)$ be an agreeable subgroup of $\GL_2(\Zhat)$. Fix a subgroup $B$ of $\calG$ such that $[\calG,\calG]\subseteq B \subseteq \calG \cap \operatorname{SL}_2(\Zhat)$. 
\begin{defn} \label{D:family}
The \defi{family} attached to the pair $(\calG,B)$, denoted $\scrF(\calG,B)$, is the set of open subgroups $G$ of $\calG$ such that $\operatorname{det}(G)=\Zhat^{\times}$ and $G\cap \operatorname{SL}_2(\Zhat) = B$. 
\end{defn}
Assume $\scrF(\calG,B)$ is nonempty and fix $G\in \scrF(\calG,B)$. In \S \ref{S:families}, we show that the family $\scrF(\calG,B)$ consists of the groups 
\begin{align}\label{E:gamma}
G_\gamma:= \{g \in \calG \colon g\cdot G = \gamma(\operatorname{det}(g))\} 
\end{align}

\noindent where $\gamma \colon \Zhat^{\times} \longrightarrow \calG/G$ a continuous homomorphism. Since the genus of an open subgroup of $\GL_2(\Zhat)$ is determined by its intersection with $\SL_2(\Zhat)$, all the groups in $\scrF(\calG,B)$ have the same genus.

\indent The families $\scrF(\calG,B)$ were first introduced in \cite{explicitopenimage}*{\S 14}. We will inspect the families of modular curves in more detail in \S \ref{S:families}. For the remainder of the paper, by a \defi{family of groups}, we mean a nonempty family attached to an arbitrary pair $(\calG,B)$.

\indent An open subgroup $G$ clearly lies in the family $\scrF(\calG, G \cap \operatorname{SL}_2(\Zhat))$, where $\calG$ is the agreeable closure of $G$. Consider the associated modular curve $X_G$ and the natural map $\pi_{\calG/G} \colon X_G \rightarrow X_\calG$. The group $\calG$ acts on $X_G$, and the restricted action of $G$ is trivial. Hence, there is an action of $\calG/G$ on $X_G$. In fact, we have $\operatorname{Aut}(X_G/X_\calG) \simeq \calG/G$, where $\operatorname{Aut}(X_G/X_\calG)$ is the group of automorphisms $f$ of the curve $X_G$ that satisfy $\pi_{\calG/G} \circ f = \pi_{\calG/G}$. The family $\scrF(\calG, G \cap \operatorname{SL}_2(\Zhat))$ consists of groups of the form $G_\gamma$, where $\gamma \colon \Zhat^{\times} \rightarrow \calG/G$ is a continuous homomorphism. Precomposing $\gamma$ with the cyclotomic character $\chi_{\operatorname{cyc}}$, we obtain a homomorphism 
\[\xi \colon \operatorname{Gal}(\QQ^{\operatorname{ab}}/\QQ) \longrightarrow \calG/G\]
which can be viewed as a $1$-cocycle of $X_G$ using the identification with curve automorphisms specified in \S \ref{S:twist}. We will show in \S \ref{S:twist} that, by twisting the curve $X_G$ with the cocycle $\xi$, we obtain a curve $(X_G)_{\xi}$ and $(X_G)_{\xi} = X_{G_\gamma}$. Hence, a family $\scrF(\calG, B)$ can be viewed as a family of twists of modular curves. 

In the rest of the paper, the terms \defi{family of groups} and \defi{family of curves} will be used interchangeably, and both will refer to a family of the form $\scrF(\calG, B)$. We prove Theorem \ref{T:ilk} as part of the following.

\begin{thm}\label{T:main families1}
Fix a non-negative integer $g$. 
\begin{enumerate}
    \item  There are only finitely many families of modular curves of genus $g$, each of which is effectively computable. (See Corollary \ref{C:twist} + Theorem \ref{T:main familiesgroups}.)
    \item Let $G$ be an open subgroup of $\GL_2(\Zhat)$ of genus $g$ with full determinant and $-I\in G$. There is an effective algorithm utilizing families of modular curves, which takes as input the group $G$ and outputs a projective curve $C_{G}\subseteq \PP_\QQ^r$ for some $r> 0$ such that $C_{G}$ is isomorphic to $X_{G}$. (See Section \ref{S:algorithm}.)

\end{enumerate}

\end{thm}

For the implementation of the algorithm of Theorem $\ref{T:main families1}$, first, we compute the finitely many families of genus $g$, choose a representative group $H\in \scrF(\calG,B)$ for each family, and precompute a projective model for $X_H$. This reduces the computation of $C_{G}$ to identifying a family $\scrF(\calG,B)$ that contains $G$, computing the continuous homomorphism $\gamma:\Zhat^{\times}\rightarrow \calG/H$ such that $G=H_\gamma$, and then twisting the projective model of $X_H$ with respect to the cocycle $\gamma\circ \chi_{cyc}$. Computationally, twisting the model reduces to a collection of linear algebra and group theory computations that can be carried out in computer algebra systems. Refer to \S \ref{S:algorithm} for a detailed exposition of our algorithm.

Theorem \ref{T:main families1} will be proved in \S \ref{S:finiteness} and \S \ref{S:twist}. Table \ref{tab:numbers} shows the number of families and agreeable subgroups for small $g$.

\begin{center}
    \begin{table}
        \centering
        \begin{tabular}{c|c|c}
      
        \toprule
        \textbf{Genus} & \textbf{ Families}  & \textbf{Agreeable Groups} \\
        \midrule
        $0$ & $638$  & $418$ \\
        $1$ & $1753$  & $1078$ \\
        $2$ & $1209$  & $885$ \\
        $3$ & $3865$  & $2244$ \\
        $4$ & $1573$ & $1151$ \\
        $5$ & $6181$ & $3659$ \\
        $\le 6$ & $15943$ & $9998$ \\
        $\le 12$ & $48819$ & $30233$ \\
        $\le 24$ & $166141$ & $95981$ \\
    
        \bottomrule
      \end{tabular}
        \caption{Number of nonempty families, indexed by pairs $(\calG,B)$ up to simultaneous conjugacy in $\GL_2(\Zhat)$, and agreeable groups up to conjugacy, for small genus.}
        \label{tab:numbers}
    \end{table}
\end{center}

\subsection{Example: A family of quadratic twists}\label{SS:example}

\indent Consider the family of modular curves $\scrF(\calG,B)$ where $\calG$ and $B$ are given as
\[
\calG:=\big\langle \begin{psmallmatrix}
    3& 10\\
    6 & 5
\end{psmallmatrix},
\begin{psmallmatrix}
    1& 6\\
    6 & 15
\end{psmallmatrix},
\begin{psmallmatrix}
    7& 2\\
    2 & 9
\end{psmallmatrix},
\begin{psmallmatrix}
    5& 0\\
    9 & 3
\end{psmallmatrix}
\big\rangle \subset \GL_2(\ZZ/16\ZZ),
\]
\[
B:=\big\langle \begin{psmallmatrix}
    7& 0\\
    0 & 7
\end{psmallmatrix}\big\rangle \subset \SL_2(\ZZ/8\ZZ).
\]
We view $\calG$ as an open subgroup of $\GL_2(\Zhat)$ by identifying it with $\pi_{16}^{-1}(\calG)$ where $\pi_{16}:\GL_2(\Zhat)\rightarrow \GL_2(\ZZ/16\ZZ)$ and similarly for $B$. The group
\[
H:=\big\langle \begin{psmallmatrix}
    9& 14\\
    6 & 7
\end{psmallmatrix},
\begin{psmallmatrix}
    7& 10\\
    5 & 9
\end{psmallmatrix},
\begin{psmallmatrix}
    5& 6\\
    7 & 3
\end{psmallmatrix}\big\rangle \subset \GL_2(\ZZ/16\ZZ).
\]

\noindent lies in the family $\scrF(\calG,B)$, has index $192$ in $\GL_2(\Zhat)$ and has genus $5$. The modular curve $X_H$ has the canonical model $C_H \subseteq \PP_\QQ^4$ given by:
\begin{align*}
    -x_1x_4 - x_2^2 + x_3^2=0\\
    -2x_1x_4 + 2x_2^2 + x_5^2=0\\
    -2x_1^2 + 2x_3x_5 + x_4^2=0.
\end{align*}

The model of $X_H$ given above is the canonical model i.e., comes from the cusp forms $S_{2,H}$ (see \S\ref{S:Background}) which is a $\QQ$-vector space of dimension 5. The agreeable closure $\calG$ has index $96$ in $\GL_2(\Zhat)$ and so $\calG/H$ is isomorphic to the cyclic group of order $2$. The nontrivial element of $\calG/H$ acts on $S_{2,H}$ (with respect to our choice of basis) via the matrix 
\[
M=\begin{psmallmatrix}
    -1& 0&0&0&0\\
    0 & 1&0&0&0\\
    0&0&-1&0&0\\
    0&0&0&-1&0\\
    0&0&0&0&-1
\end{psmallmatrix}
\]

Let $\gamma:\Zhat^{\times}\rightarrow\calG/H$ be a continuous homomorphism and let $H_\gamma$ be the corresponding group in $\scrF(\calG,B)$. Since $\calG/H$ is cyclic of order two, there exists a squarefree integer $d$ such that $\chi_{\operatorname{cyc}}(\Gal_{K_\gamma})=\ker \gamma$, where $K_\gamma:=\QQ(\sqrt{d})$. The curves $X_H$ and $X_{H_\gamma}$ are isomorphic over $K_\gamma$. For nontrivial $\gamma$, the matrix cocycle defining the twist $C_{H_\gamma}$ of the model $C_H$ is given by the map $\bar\xi:\Gal(K_\gamma/\QQ)\rightarrow \GL_5(K_\gamma)$, $ \sigma\ne 1\mapsto M$. For trivial $\gamma$, we take $d=1$ and the cocycle is trivial. A computation shows that the curve given by the equations
\begin{align*}
    -x_1x_4 - dx_2^2 + x_3^2=0\\
    -2x_1x_4 + 2dx_2^2 + x_5^2=0\\
    -2x_1^2 + 2x_3x_5 + x_4^2=0.
\end{align*}
is isomorphic to $X_{H_\gamma}$. Hence, the set of equations given above as we vary the squarefree integer $d$, corresponds to the family of modular curves $\scrF(\calG,B)$. As a particular example, consider the following group.
\[
G=\big\langle \begin{psmallmatrix}
    429& 214\\
    270 & 211
\end{psmallmatrix},
\begin{psmallmatrix}
    633& 548\\
    824 & 505
\end{psmallmatrix},
\begin{psmallmatrix}
    425& 275\\
    630 & 871
\end{psmallmatrix},
\begin{psmallmatrix}
    663& 909\\
    610 & 913
\end{psmallmatrix}\big\rangle \subset \GL_2(\ZZ/944\ZZ).
\]

The curves $X_H$ and $X_G$ are isomorphic over the number field $K:=\QQ(\sqrt{118})\subset \QQ(\zeta_{944})$. Using our implementation, we find the curve $C\subseteq \PP_\QQ^4$ given by the equations
\begin{align*}
    -x_1x_4 - 118x_2^2 + x_3^2=0\\
    -2x_1x_4 + 236x_2^2 + x_5^2=0\\
    -2x_1^2 + 2x_3x_5 + x_4^2=0.
\end{align*} 
is isomorphic to $X_G$. 

On our machine, the precomputation of the family $\scrF(\calG,B)$, the representative $X_H$, and the model $C_H$ took $0.50$ seconds using Zywina's method \cite{explicitopenimage}. After this precomputation, our implementation took 0.41 seconds to compute the model $C$ of $X_G$ and the associated $j$-map $X_G \rightarrow \PP_\QQ^1$. On the same machine, Zywina’s implementation took $22.79$ seconds to compute a model for the modular curve $X_G$ (without the $j$-map). For modular curves of higher level, the computation of $j$-maps is particularly time consuming, as it involves finding relations between the $j$-invariant and cusp forms whose coefficients lie in $\QQ(\zeta_N)$, where $N$ is the level of $G$. We avoid this step via twisting in our algorithm. As a result, our algorithm remains efficient even as the level of the input curve $X_G$ increases. Although we do not compute $q$-expansions to find the models, it should be noted that the current implementation of our algorithm can still be used to compute certain subspaces of $M_{k,G}$ (for instance, $S_{2,G}$ in the above example) via twisting.

\subsection{An example of cubic twists} 
\indent Consider the following group.
\[
G:=\big\langle \begin{psmallmatrix}
    41 & 160\\
    4 & 61
\end{psmallmatrix},
\begin{psmallmatrix}
    155& 77\\
    135 & 60
\end{psmallmatrix},
\begin{psmallmatrix}
    121& 93\\
    13 & 148
\end{psmallmatrix}
\big\rangle \subset \GL_2(\ZZ/182\ZZ).
\]

\noindent We find that $G$ is conjugate to a group that lies in the family $\scrF(\calG,B)$ where $\calG$ and $B$ are given as
\[
\calG:=\big\langle \begin{psmallmatrix}
    8& 3\\
    13 & 3
\end{psmallmatrix},
\begin{psmallmatrix}
    23& 14\\
    8 & 1
\end{psmallmatrix},
\begin{psmallmatrix}
    11& 13\\
    17 & 16
\end{psmallmatrix}
\big\rangle \subset \GL_2(\ZZ/26\ZZ),
\]
\[
B:=\big\langle \begin{psmallmatrix}
    17& 22\\
    18 & 5
\end{psmallmatrix},
\begin{psmallmatrix}
    23& 16\\
    14 & 3
\end{psmallmatrix}\big\rangle \subset \SL_2(\ZZ/26\ZZ).
\]
We have a representative $H$ in the family $\scrF(\calG,B)$ given by
\[
H:=\big\langle \begin{psmallmatrix}
    13& 4\\
    22 & 23
\end{psmallmatrix},
\begin{psmallmatrix}
    15& 18\\
    6 & 25
\end{psmallmatrix}\big\rangle \subset \GL_2(\ZZ/26\ZZ).
\]

\noindent The groups $G$ and $H$ have index $84$ in $\GL_2(\Zhat)$ and have genus $5$. The modular curve $X_H$ has the canonical model $C_H \subseteq \PP_\QQ^4$ defined by the polynomials 
\begin{align*}
    f_1 &= x_2x_4 + x_3x_5, \\[6pt]
f_2 &= -x_1^2 + x_1x_2 + x_1x_3 - x_1x_4 + x_1x_5 - x_2x_3 + x_2x_4 + x_4x_5, \\[6pt]
f_3 &= -x_1^2 + x_2^2 - x_2x_3 - x_2x_4 + x_2x_5 + x_3^2 - x_3x_4 + x_4^2 + x_4x_5 + x_5^2.
\end{align*}
Using our implementation, we find the curve $C\subseteq \PP_\QQ^4$ defined by the polynomials
\begin{align*}
f_1' = & -4x_1^2 + 14x_1x_2 + 6x_1x_3 - 26x_1x_4 + 26x_1x_5 
       + 17x_2^2 - 42x_2x_3 \\
       &+ 5x_2x_4 - 50x_2x_5 
       + 27x_3^2 + 24x_3x_4 + 21x_3x_5 
       + 5x_4^2 + 26x_4x_5 + 5x_5^2, \\[3pt]
f_2' = & 4x_1^2 - 14x_1x_2 + 24x_1x_3 + 5x_1x_4 + 7x_1x_5 
       + 10x_2^2 + 3x_2x_3 \\
       &- 11x_2x_4 + 8x_2x_5 
       + 6x_3x_4 - 24x_3x_5 + x_4^2 - 17x_4x_5 - 11x_5^2, \\[3pt]
f_3' = & 6x_1^2 - 21x_1x_2 + 29x_1x_3 + 4x_1x_4 + 7x_1x_5 
       + 15x_2^2 + x_2x_3\\
       &- 13x_2x_4 + 5x_2x_5 
       + 2x_3x_4 - 29x_3x_5 - 2x_4^2 - 15x_4x_5 - 13x_5^2.
\end{align*}
is isomorphic to $X_G$. The curves $X_H$ and $X_G$ are isomorphic over the number field $K:=\QQ(z)\subset \QQ(\zeta_{182})$ with the defining polynomial $f(x)=x^3-5x^2+6x-1$. The model $C_H$ is obtained from the cusp forms $S_{2,H}$ (see \S \ref{S:Background}) which is a $\QQ$-vector space of dimension $5$. The agreeable closure $\calG$ has index $28$ in $\GL_2(\Zhat)$ and $\calG/H$ is isomorphic to the cyclic group of order $3$. Choose the generator $\sigma\in\Gal(K/\QQ)$ with $\sigma(z)=z^2-4z+4$, and let $g\in \calG/H$ be its image under $\xi$ (see \S \ref{SS:famsfirst}). The generator $g$ of $\calG/H$ acts on $S_{2,H}$ (with respect to our choice of basis and generator) via the matrix 
\[
M=\begin{psmallmatrix}
    1& 0&0&0&0\\
    0 & 0&0&0&1\\
    1&0&-1&1&0\\
    0&0&-1&0&0\\
    1&-1&0&0&-1
\end{psmallmatrix}
\]

The associated matrix cocycle (which is a homomorphism) is given by the map $\bar\xi:\Gal(K/\QQ)\rightarrow \GL_5(K)$, $\sigma\mapsto M$. Hilbert 90 implies that $\bar\xi$ is a coboundary, and the matrix 
\[
A:=\begin{psmallmatrix}
1 & -1 & 0 & 1 & -1 \\[2pt]
\tfrac{2z - 1}{7} & \tfrac{3z^2 - 8z - 2}{7} & 0 & \tfrac{3z^2 - 13z + 11}{7} & \tfrac{-2z + 1}{7} \\[2pt]
\tfrac{2z - 1}{7} & \tfrac{-2z + 1}{7} & \tfrac{-3z^2 + 15z - 12}{7} & \tfrac{2z - 1}{7} & \tfrac{3z^2 - 8z - 2}{7} \\[2pt]
\tfrac{2z^2 - 6z - 1}{7} & \tfrac{-2z^2 + 6z + 1}{7} & \tfrac{3z^2 - 6z - 3}{7} & \tfrac{2z^2 - 6z - 1}{7} & \tfrac{4z^2 - 15z + 10}{7} \\[2pt]
\tfrac{2z^2 - 8z + 7}{7} & \tfrac{z^2 - 7z + 5}{7} & 0 & \tfrac{-4z^2 + 13z - 2}{7} & \tfrac{-2z^2 + 8z - 7}{7}
\end{psmallmatrix}
\]
satisfies $\bar\xi(\tau) = \tau(A)A^{-1}$ for all $\tau \in \Gal(K/\QQ)$. Hence, $A$ encodes how to pass between two $\QQ$-structures defined on $S_{2,H}\otimes_\QQ K$ corresponding to $S_{2,H}$ and $S_{2,G}$ allowing us to twist the curve $C_H$ and obtain the curve $C$.

\subsection{Motivation} In \cite{Mazur72}, Mazur describes the following program which serves as motivation for computing projective models of modular curves and for classification problems related to modular curves:
\begin{programb}
 Given a number field $K$ and a subgroup $H$ of $\operatorname{GL}_2(\Zhat)=\prod_p \operatorname{GL}_2(\ZZ_p)$ classify all elliptic curves $E/K$ whose associated Galois representation on torsion points maps $\operatorname{Gal}(\overline{K}/K)$ \underline{into} $H\subseteq \operatorname{GL}_2(\Zhat)$.
\end{programb}

Let $G_1$ and $G_2$ be the groups $\pm\rho_{E_1}^{\ast}(\Gal_\QQ)$ and $\pm\rho_{E_2}^{\ast}(\Gal_\QQ)$, where $E_1$ and $E_2$ are elliptic curves with j-invariants $-7\cdot 11^3$ and $-7 \cdot 137^3 \cdot 2083^3$, respectively. Note that these groups are well-defined up to conjugacy in $\GL_2(\Zhat)$.

In \cite{explicitopenimage}, it is conjectured that if $G \subset \GL_2(\Zhat)$ is an open subgroup with surjective determinant containing $-I$, and if $X_G$ has genus at least $54$, and $G$ is not conjugate to $G_1$ or $G_2$ in $\GL_2(\Zhat)$, then $X_G$ contains no non-CM rational points over $\QQ$. Hence, conjecturally, explicitly computing the families of modular curves up to genus $53$, along with extending our algorithm to such families, allows us to compute a projective model for any modular curve over rationals that contains a non-CM rational point (except $X_{G_1}$ and $X_{G_2}$ whose rational points are understood).

Following the conjectures of Zywina, we suggest the following challenging program:
\begin{program} One can consider the following steps to resolve Mazur's Program B:
    \begin{enumerate}
        \item Prove Serre's uniformity problem.
        \item Classify all rational points on a finite number of special modular curves as described in Section 14 of \cite{explicitopenimage}.
        \item Classify all congruence subgroups of $\SL_2(\ZZ)$ (in the sense of \cite{MR2016709}) up to genus 53 (or genus $\beta$ as in Lemma 14.7 in \cite{explicitopenimage}).
        \item Compute all families of modular curves up to the genus mentioned in (3), in the sense of \S \ref{S:families}.
        \item Investigate the behavior of rational points on the mentioned families. 
    \end{enumerate}
\end{program}
Note that there has been much progress towards proving Serre's uniformity problem in which Serre asks whether for all primes $l>37$, the mod $l$ representation $\rho_{E,l}$ is surjective or not. If the image is not surjective, then $E$ gives rise to a non-CM rational point on the modular curve $X_G$, where $G$ is a maximal subgroup of $\GL_2(\ZZ/l\ZZ)$. Mazur \cite{ratiso,ikincimazur} completely described the cases where $G$ is the Borel subgroup or one of the exceptional subgroups of $\GL_2(\ZZ/l\ZZ)$. Bilu, Parent and Rebolledo \cite{BPR} showed that when $G$ is the normalizer of split Cartan subgroup, then $X_G$ has no non-CM rational points. The only remaining case is the normalizer of non-split Cartan subgroups and the associated modular curves $X_{ns}^+(l)$. Using Chabauty methods, Balakrishnan et al. \cite{cursed,MR4589060} determined rational points on $X_{ns}^+(l)$ for some small primes $l$. 

\subsection{Related results}
\indent There is a lot of work on modular curves and Galois representations attached to elliptic curves over $\QQ$. Here, we mention some recent related results.\\
\indent In \cite{explicitopenimage}, Zywina describes a practical algorithm that computes the image of $\rho_E$ up to conjugacy in $\operatorname{GL}_2(\Zhat)$. Assuming some conjectures, they also give a complete classification of the groups $\rho_E(\operatorname{Gal}(\overline{\QQ}/\QQ)) \cap \operatorname{SL}_2(\Zhat)$. The methods developed by Zywina to achieve this result include an algorithm to compute models of modular curves $X_G$. They also introduce the notion of a family of modular curves and interpret their results in this language. The notions and methods introduced by Zywina form an important basis for our work. In \cite{openimagenumberfields}, they describe an analogous algorithm that works for elliptic curves over number fields. \\
\indent Rakvi \cite{rakvi} has given a classification of genus 0 modular curves $X_G$ over $\QQ$ such that $X_G \cong \PP_{\QQ}^1$, in terms of families of abelian twists.\\
\indent In \cite{MR3671434}, the authors determine all open subgroups $G$ of prime power level for which $X_G(\QQ)$ is infinite. This work also provides a classification for possible images of $l$-adic Galois representations arising from elliptic curves for almost all $j$-invariants. \\
\indent In \cite{rousezureick}, Rouse and Zureick-Brown give a classification of the possible $2$-adic images of Galois representations associated to non-CM elliptic curves over $\QQ$.
Rouse, Sutherland and Zureick-Brown \cite{rszb} investigate the $l$-adic images for $l=3,5,7,11$. Building on their work, Balakrishnan et al.\ \cite{balakrishnan2025} completed the classification of $3$-adic images.
In \cite{furio7adic}, Furio and Lombardo reduce the $7$-adic classification to computing the rational points on a single plane quartic.
In forthcoming work, we compute the rational points on the plane quartic of \cite{furio7adic} and on $X_{ns}^+(25)$ using quadratic Chabauty.
Arango-Pi\~neros and Zureick-Brown \cite{azb7adic} determine the rational points on the two remaining modular curves of level $49$ and genus $9$, completing the classification of $7$-adic images for non-CM elliptic curves over $\QQ$.
Arango-Pi\~neros, Najman and Zureick-Brown \cite{anzb11adic} announce that $X_{ns}^+(121)$ has exactly seven rational points, all CM, which would complete the corresponding $11$-adic classification.
In recent preprints, Jaswal \cites{jaswal25,jaswal121} reports nine rational points on $X_{ns}^+(25)$ and seven on $X_{ns}^+(121)$ recovering all known CM points, and Nguyen \cite{nguyen7adic} reports four on the plane quartic of \cite{furio7adic}.
For a prime $p>7$ and a non-CM elliptic curve $E/\QQ$ whose mod $p$ image lies in the normalizer of a non-split Cartan subgroup, Bisatt, Furio and Lombardo \cite{lombardo} show that the $p$-adic image is the full inverse image of such a normalizer modulo $p^n$ for some $n\geq 1$.

Mayle and Rouse \cite{jeremy} give an algorithm to determine $X_G(\QQ)$ from an explicit nonconstant map over $\QQ$ to an elliptic curve of rank zero.

\subsection{Implementation}
The implementation of our algorithm can be found in the repository \cite{repo}. 

\begin{center}
    \href{https://github.com/eekarabiyik/twist}{https://github.com/eekarabiyik/twist}
\end{center}
All computations are done in \texttt{Magma} \cite{Magma}. This repository also includes all families of modular curves up to genus $12$, relevant representative groups, and models. There are hundreds of thousands of such objects and we are not able to express them in tables. 

\indent In Section \ref{S:algorithm}, we describe our algorithm including the precomputation part. While the precomputation part is computationally intense, it is only a one time computation and the rest of our algorithm is efficient.

\subsection{Notation}
$\Zhat$ is the profinite ring obtained by taking the inverse limit of $\ZZ/N\ZZ$ over all $N\in \NN$. Similarly, for $M\in\NN$ with $M>1$, $\ZZ_M$ is the profinite ring obtained by taking the inverse limit of $\ZZ/M^s\ZZ$ where $s$ ranges over $\NN$. There are natural isomorphisms 
\[
\ZZ_M \cong \prod_{l \vert M} \ZZ_l \quad \text{and} \quad  \Zhat \cong \prod_{l}\ZZ_l
\]
where the product runs over the prime numbers $l$. The reduction modulo $N$ homomorphism $\Zhat \rightarrow \ZZ/N\ZZ$  induces homomorphism $\pi_N:\operatorname{GL}_2(\Zhat) \rightarrow \operatorname{GL}_2(\ZZ/N\ZZ)$. The level of an open subgroup $G$ of $\operatorname{GL}_2(\Zhat)$ is the smallest positive integer $N$ such that $G$ is the inverse image of the reduction modulo $N$ map $\pi_N:\operatorname{GL}_2(\Zhat) \rightarrow \operatorname{GL}_2(\ZZ/N\ZZ)$ i.e., $G=\pi_N^{-1}(\pi_N(G))$. In this case, $G$ is uniquely determined by $\pi_N(G)$, we will often denote it by $G$ as well, when the level is understood. Similarly, the level of an open subgroup of $\operatorname{GL}_2(\ZZ_M)$ is the smallest positive integer $N$ that divides a power of $M$ and such that $G$ is equal to the inverse image of its image under the reduction modulo $N$ map $\operatorname{GL}_2(\ZZ_M) \rightarrow \operatorname{GL}_2(\ZZ/N\ZZ)$. The levels of an open subgroup of $\operatorname{SL}_2(\Zhat)$ and $\operatorname{SL}_2(\ZZ_M)$ are defined similarly.

\indent For $1< N\in \NN$, we let $G_N$ be the image of $G$ under the homomorphism $\GL_2(\Zhat)\to \GL_2(\ZZ_N)$ arising from the natural projection map. If $G$ has level $N>1$, then $G=G_N \times \prod_{\ell \nmid N} \GL_2(\ZZ_\ell)$. We denote by $G(N)$ the image of $G$ under the homomorphism $\GL_2(\Zhat) \rightarrow \GL_2(\ZZ/N\ZZ)$.\\
\indent Unless stated otherwise, open subgroups $G$ of $\operatorname{GL}_2(\Zhat)$ are assumed to satisfy $\det (G) = \Zhat^{\times}$ and $-I \in G$. 
 \subsection{Overview of the Paper}
 
\indent In \S \ref{S:Background}, we review the background material on modular curves and modular forms. In \S \ref{S:agreeable}, we will discuss agreeable subgroups, define the agreeable closure of a subgroup $G$ and discuss how to compute it. In \S \ref{S:families}, we define a family of groups associated to a pair $(\calG,B)$, denoted by $\scrF(\calG,B)$. In \S \ref{S:finiteness}, we prove the finiteness of agreeable subgroups of a fixed genus and use this to deduce our main theorem, that modular curves over $\QQ$ lie in finitely many families of abelian twists. In \S \ref{S:twist}, we will describe the cocycles arising from families $\scrF(\calG,B)$, show that $\scrF(\calG,B)$ is a family of twists of modular curve and describe how to twist projective models of modular curves. In \S \ref{S:algorithm}, we finally describe an algorithm for computing a model of any modular curve $X_G$ of fixed genus $g$. In the appendix, we assume $-I\notin G$ and explain how to compute quotients of modular forms, certify them using explicit precision bounds, and construct a model for the generic fiber of the universal elliptic curve associated to $G$. 
 
 \subsection*{AI disclosure} The main text of this paper forms the first part of the author’s PhD thesis, which was submitted on 30 April 2026. We used an AI assistant only for correcting typographical errors and improving the clarity of the prose. During the typographical check, the AI assistant pointed out an error in the order of multiplication of Hilbert 90 matrices in \S \ref{SS:getting a basis}. This error was fixed by the author independently.

 \section{Modular Forms and Modular Curves} \label{S:Background}
 \indent The goal of this section is to state some well known facts about the theory of modular forms and introduce modular curves. We will mostly use the language of \cite{explicitopenimage}. For the rest of the section, let $G \subseteq \operatorname{GL}_2(\Zhat)$ be an open subgroup. For any $N$ divisible by the level of $G$, the projection $\pi_N: \operatorname{GL}_2(\Zhat) \rightarrow \operatorname{GL}_2(\ZZ/N\ZZ)$ gives a group whose inverse image is the open subgroup $G$ i.e., $G=\pi_N^{-1}(\pi_N(G))$. We will often abuse notation and denote by $G$ both the open subgroup of $\GL_2(\Zhat)$ and its image under $\pi_N$. Considering $G$ as a subgroup of $\operatorname{GL}_2(\ZZ/N\ZZ)$, we let $\Gamma_G$ be the subgroup of $\operatorname{SL}_2(\ZZ)$ consisting of matrices that are modulo $N$ congruent to an element of $G \cap \operatorname{SL}_2(\ZZ/N\ZZ)$.

 \subsection{Setting the stage}
  
\indent Let $\HH$ be the complex upper plane. The group $\operatorname{SL}_2(\ZZ)$ acts on the extended complex upper half-plane $\HH^{\ast}= \HH \cup \QQ \cup \{\infty\} $ by linear fractional transformations. Fix a congruence subgroup $\Gamma$ of $\SL_2(\ZZ)$ such that $-I \in \Gamma$. For a positive integer $N$, define the primitive $N$-th root of unity $\zeta_N:=e^{2\pi i/N}$ in $\CC$.  

The quotient $\calX_\Gamma:=\Gamma\backslash \HH^*$ is a smooth compact Riemann surface \cite{MR2112196}*{\S 2}. Let $g$ be the genus of the Riemann surface $\calX_\Gamma$. Let $w$ be the width of the cusp $\infty$, i.e. the smallest positive integer such that $\left[\begin{smallmatrix}
1 & w \\0 &1
\end{smallmatrix}\right] \in \Gamma$. Let $q_w:=e^{2\pi i \tau /w}$.

Let $k\geq 0$ be a natural number. For a meromorphic function $f$ on $\HH$ and a matrix $\gamma\in \GL_2(\RR)$ with positive determinant, we define the \emph{slash operator} of weight $k$ on $f$ by $(f|_k\gamma)(\tau):=\det(\gamma)^{k/2}(c\tau+d)^{-k}f(\gamma\tau)$. Let $P_1,\ldots, P_r$ be the cusps of $\calX_\Gamma$. Let $Q_1,\ldots, Q_s$ be the elliptic points of $\calX_\Gamma$ and denote their orders by $e_1,\ldots, e_s$, respectively. Each $e_i$ is either $2$ or $3$.  Let $v_2$ and $v_3$ be the number of elliptic points of $\calX_\Gamma$ of order $2$ and $3$, respectively.

\subsection{Modular forms}\label{SS:modularformsgeneral}

A modular form of weight $k\geq 0$ with respect to $\Gamma$ is a holomorphic function of $\HH$ such that $f|_k\gamma=f$ for all $\gamma \in \Gamma$, and at the cusps it satisfies the usual growth condition. We denote the set of modular forms with respect to $\Gamma$ by $M_k(\Gamma)$. It is a finite dimensional complex vector space. Let $f\in M_k(\Gamma)$ be a modular form. We have a unique $q$-expansion of $f$ (at the cusp $\infty$) given by 
\[
f(\tau) = \sum_{n=0}^\infty a_n(f)\, q_w^{n}
\]
where $a_n(f)\in \CC$. The spaces of modular forms of different weights $k$ form a graded $\CC$-algebra denoted by 
\[
R_{\Gamma}: = \bigoplus_{k \geq 0} M_{k}(\Gamma).
\]  
 The graded ring $R_\Gamma$ is finitely generated as a $\CC$-algebra. One can focus on modular forms $f$ whose $q$-expansion (at the cusp $\infty$) has coefficients in a certain subring $S$ of $\CC$, which we denote by  $M_k(\Gamma,S)$. It has a natural structure as an $S$-module.\\
\indent In particular, consider the principal congruence subgroups $\Gamma(N)$ given by
\[
\Gamma(N)\colonequals \bigg\{ \begin{bmatrix}
a & b\\
c & d 
\end{bmatrix} \in \operatorname{SL}_2(\ZZ) \quad \bigg| \quad a,d \equiv 1 \;(\operatorname{mod } N)\;\; b,c \equiv\; 0 \;(\operatorname{mod } N)\bigg\}.
\]
The $\QQ(\zeta_N)$-vector space $M_{k}(\Gamma(N),\QQ(\zeta_N))$ will be of special interest to us.

\subsection{Modular forms and differential forms}
\label{SS:global sections}

Fix an even integer $k\geq 0$. Take any modular form $f\in M_k(\Gamma)$. By definition, $f$ satisfies $f(\gamma\tau)=(c\tau+d)^{k}f(\tau)$ for all $\gamma \in \Gamma$, which is equivalent to $f(\gamma \tau) d( \gamma \tau)^{k/2} = f(\tau) (d\tau)^{k/2}$. Therefore, $f$ gives a differential form 
\begin{align} \label{E:form on H}
 (2\pi i)^{k/2}f(\tau) \, (d\tau)^{k/2} = w^{k/2} \bigg(\sum_{n=0}^\infty a_n(f)\, q_w^{n}\bigg) \bigg(\frac{dq_w}{q_w} \bigg)^{k/2}
\end{align}
on $\HH$ and this induces a meromorphic differential $k/2$-form $\omega_f$, associated to $f$, on $\calX_\Gamma$. For details, see \cite{MR2112196}*{\S 3.3}.

Let $\calD_k$ be the divisor 
\begin{align}\label{E:divisor}
\sum_{i=1}^{r}k/2 \cdot P_i + \sum_{i=1}^s \lfloor k/2\cdot(1-1/e_i)\rfloor\cdot Q_i 
\end{align}

of $\calX_\Gamma$; it is supported on cusps and elliptic points \cite{explicitopenimage}*{\S 4.4}. For any modular form $f\in  M_k(\Gamma)$, we have $\operatorname{div}(\omega_f) + \calD_k \geq 0$. Writing $\Omega$ for the sheaf of differentials on the curve under consideration, we obtain a map of complex vector spaces
\[
\psi_k\colon M_k(\Gamma) \rightarrow H^0(\calX_\Gamma,\Omega^{\otimes k/2}(\calD_k))
\] 
sending $f$ to $\omega_f$.

Moreover, any differential form in $H^0(\calX_\Gamma,\Omega^{\otimes k/2}(\calD_k))$ pulls back to a differential form $(2\pi i)^{k/2}f(\tau)(d\tau)^{k/2}$ on $\HH$ as in (\ref{E:form on H}). The form $f$ is holomorphic and satisfies the growth conditions at cusps (due to being in the space $H^0(\calX_\Gamma,\Omega^{\otimes k/2}(\calD_k))$). Therefore, the map $\psi_k$ is an isomorphism of vector spaces.

The groups $\Gamma$ that we consider contain $-I$, which means that $M_k(\Gamma)=0$ when $k$ is odd. Combining the isomorphisms $\psi_k$ for even $k\geq0$, we get an isomorphism of graded $\CC$-algebras:
\[
\psi\colon R_{\Gamma} \xrightarrow{\sim} \bigoplus_{\substack{k \geq 0\\ k\ \mathrm{even}}} H^0(\calX_{\Gamma},\Omega^{\otimes k/2}(\calD_k)).
\]



\subsection{Actions}
\indent Fix positive integers $k> 1$ and $N$. The principal congruence subgroup $\Gamma(N)$ is normal in $\operatorname{SL}_2(\ZZ)$ and the weight-$k$ operator gives a right action of $\SL_2(\ZZ/N\ZZ)$ on $M_k(\Gamma(N))$. Let $f=\sum_{n=0}^\infty a_n(f) q_N^n$ be a modular form in $M_k(\Gamma(N))$. Let $\sigma$ be a field automorphism of $\CC$, it acts on the coefficients of $f$ and gives rise to a unique weight-$k$ modular form $\sigma(f)$, i.e. the $q$-expansion of $\sigma(f)$ is given by $\sum_{n=0}^\infty \sigma(a_n(f)) q_N^n$.

Consider the isomorphism $(\ZZ/N\ZZ)^\times \xrightarrow{\sim} \Gal(\QQ(\zeta_N)/\QQ)$, $d\mapsto \sigma_d$, where $\sigma_d(\zeta_N)=\zeta_N^d$. The following lemma gives an action of $\GL_2(\ZZ/N\ZZ)$ on $M_{k}(\Gamma(N),\QQ(\zeta_N))$ viewed as a $\QQ$-vector space.   

\begin{lemma} \label{L:star action}
There is a unique right action $*$ of $\GL_2(\ZZ/N\ZZ)$ on $M_k(\Gamma(N),\QQ(\zeta_N))$ such that the following hold:
\begin{itemize}
\item
if $A \in \SL_2(\ZZ/N\ZZ)$, then  $f*A = f|_k \gamma$, where $\gamma$ is any matrix in $\SL_2(\ZZ)$ that is congruent to $A$ modulo $N$,
\item
if $A=\left(\begin{smallmatrix}1 & 0 \\0 & d\end{smallmatrix}\right)$, then $f*A=\sigma_d(f)$.
\end{itemize}
\end{lemma}
\begin{proof}
See \cite[\S3]{BN2019}.
\end{proof}

 \indent Combining the action of Lemma $\ref{L:star action}$ for all $k$, we get a right action $*$ of $\GL_2(\ZZ/N\ZZ)$ on the graded $\QQ$-algebra $\bigoplus_{k\geq 0} M_k(\Gamma(N),\QQ(\zeta_N))$.

\subsection{The spaces $M_{k,G}$} \label{SS:MkG}
Fix a positive integer $N$. Let $G$ be a subgroup of $\GL_2(\ZZ/N\ZZ)$ such that $\det(G)=(\ZZ/N\ZZ)^{\times}$. Define the $\QQ$-vector space
\[
M_{k,G}:=M_k(\Gamma(N),\QQ(\zeta_N))^G,
\]
i.e. the subspace fixed by $G$ under the action $*$ from Lemma~\ref{L:star action}. Note that
\[
M_{k,G} \subseteq M_k(\Gamma(N),\QQ(\zeta_N))^{G\cap \operatorname{SL}_2(\ZZ/N\ZZ)} = M_k(\Gamma_G,\QQ(\zeta_N)).
\]

\indent Tensoring $M_{k,G}$ with $\QQ(\zeta_N)$ and $\CC$ give natural isomorphisms.
\begin{lemma} \label{L:tensor agreement}
The natural homomorphisms
\[
M_{k,G} \otimes_{\QQ} \QQ(\zeta_N) \to M_k(\Gamma_G,\QQ(\zeta_N)) \quad \text{ and } \quad 
M_{k,G} \otimes_{\QQ} \CC \to  M_k(\Gamma_G)
\]
are isomorphisms for $k\neq 1$.
\end{lemma}
\begin{proof}
See Lemma 4.5 in \cite{explicitopenimage}.
\end{proof}

%
%
%
%
%
%

 \subsection{Modular curves}
 Let $G$ be a subgroup of $\operatorname{GL}_2(\ZZ/N\ZZ)$ that satisfies $\operatorname{det}(G)=(\ZZ/N\ZZ)^{\times}$ and $-I \in G$. We have defined the $\QQ$-vector spaces
\[
M_{k,G}:=M_k(\Gamma(N),\QQ(\zeta_N))^G.
\]

Note that when $-I\in G$, $M_{k,G}$ is trivial for odd integers $k$. Consider the graded $\QQ$-algebra $\bigoplus _{k=0}^{\infty}M_{k,G}$.

 \begin{defn}\label{D:modcurve}
 The \defi{modular curve} $X_G$ associated to group $G$ is the $\QQ$-scheme $\operatorname{Proj}(\bigoplus _{k=0}^{\infty}M_{k,G})$. For open subgroups $G\subseteq \operatorname{GL}_2(\Zhat)$, the modular curve $X_G$ is defined using the image of $G$ under the natural projection $\pi \colon \operatorname{GL}_2(\Zhat) \rightarrow \operatorname{GL}_2(\ZZ/N\ZZ)$ where $N$ is divisible by the level of $G$. This definition is independent of the choice of $N$.
 \end{defn}
 
 \begin{remark}
The modular curve $X_G$ is a smooth, projective and geometrically irreducible curve over $\QQ$. 
 \end{remark}
 
\indent When $G=\operatorname{GL}_2(\Zhat)$, we have $X_G=\operatorname{Proj}(\bigoplus _{k=0}^{\infty}M_{k,G})=\operatorname{Proj}(\QQ[E_4,E_6])$. Its function field is $\QQ(j)$, where $j=\frac{1728E_4^3}{E_4^3-E_6^2}$ is the modular $j$-invariant. The first few terms of the $q$-expansion of $j$ are given by
\[
j(\tau)=q^{-1}+744+196884q+21493760q^2+864299970q^3+\cdots
\]
We can identify $X_{\GL_2(\Zhat)}$ with $\PP^1_{\QQ}$ coordinate $j$. It is commonly called as the \defi{\emph{j}-line}.\\
\indent Let $G\subseteq G' \subseteq \operatorname{GL}_2(\Zhat)$ be open subgroups. There is an inclusion $M_{k,G'}\subseteq M_{k,G}$, and so there is an induced map of curves $X_G\rightarrow X_{G'}$. When $G'=\operatorname{GL}_2(\Zhat)$, we obtain the \defi{absolute $j$-map} $\pi \colon X_G \rightarrow \PP^1_{\QQ}$.
 
 \subsubsection{Compatibility with other definitions}
 
\indent Consider the space of modular forms $M_k(\Gamma_G)$ as in \S \ref{SS:modularformsgeneral}, for even $k\geq0$. It is a finite dimensional complex vector space. The map $\psi_k$ identifies it with $H^0(\calX_{\Gamma_G},\calL_k)$, where $\calL_k:=\Omega^{\otimes k/2}(\calD_k)$. Consider the graded ring $R_{\Gamma_G}:=\bigoplus_{k\geq0}M_k(\Gamma_G)$ and the graded $\QQ$-algebra $R:=\bigoplus_{k\geq0}M_{k,G}$. Lemma \ref{L:tensor agreement} gives $R\otimes_{\QQ}\CC\simeq R_{\Gamma_G}$.

\indent Let $X_G$ for this paragraph be defined as in \cite{explicitopenimage}*{\S 3.3}, by explicitly describing its function field. The divisor $\calD_k$ is defined over $\QQ$, so we can view it as a divisor on $X_G$. Define the invertible sheaf $\scrL_k:=\Omega_{X_G/\QQ}^{\otimes k/2}(\calD_k)$ on $X_G$, which gives rise to $\calL_k$ on $X_G(\CC)=\calX_{\Gamma_G}$. Base change gives
\[
H^0(X_G,\scrL_k)\otimes_{\QQ}\CC
\simeq H^0(\calX_{\Gamma_G},\calL_k)\simeq M_k(\Gamma_G).
\]
The normalized map $\psi_k$ induces an isomorphism \cite{explicitopenimage}*{\S 5.1}
\[
M_{k,G} \xrightarrow[]{\sim} \operatorname{H}^0(X_G,\scrL_k ).
\]
Since $\calD_{12n}=n\calD_{12}$, the twelfth Veronese subring of $R$ is
\[
R^{(12)}:=\bigoplus_{n\geq0}M_{12n,G}
\simeq\bigoplus_{n\geq0}H^0(X_G,\scrL_{12}^{\otimes n}).
\]
The degree of $\scrL_{12}$ is $[\SL_2(\ZZ):\Gamma_G]>0$, so $\scrL_{12}$ is ample. It follows that $\operatorname{Proj}(R)\simeq\operatorname{Proj}(R^{(12)})\simeq X_G$. Hence, our definition is compatible with Zywina's definition. In particular, we identify $X_G(\CC)$ with $\calX_{\Gamma_G}$, and the base change of $\pi\colon X_G\rightarrow\PP^1_{\QQ}$ is the complex projection $\calX_{\Gamma_G}\rightarrow\PP^1_{\CC}$.

 \subsection{Models of modular curves}\label{S:Models}
 
 In this section, we briefly describe how to compute projective models of modular curves, as explained in \cite{explicitopenimage}*{\S 5}. This method uses explicit $q$-expansions of modular forms in $M_{k,G}$.

Let $G$ be as above, and let $N$ be a positive integer divisible by the level of $G$. Let $\calX_{\Gamma_G}$ be the Riemann surface associated to $\Gamma_G$ which we identify with $X_G(\CC)$. Let $P_1,\cdots , P_r$ be the cusps of $X_G(\CC)$, which are defined over $\QQ(\zeta_N)$. Let $E=\sum_{i=1}^{r}e_iP_i$ be a divisor on $X_G$ defined over $\QQ$ with $e_i\geq 0$. Let $g$ be the genus of the curve $X_G$. For a modular form $f$, let $\nu_{P_i}(f)$ be the order of its expansion $f|_k A_i$ in $q^{1/w_i}$, where $A_i\in\SL_2(\ZZ)$ sends $\infty$ to a representative of $P_i$ and $w_i$ is the cusp width. Since $-I\in\Gamma_G$, this is the order in a local parameter on $X_G$.
Define
\[
V:= \{ f\in M_{k,G} :  \nu_{P_i}(f)\geq e_i \text{ for all $1\leq i \leq r$}\}.
\]

Assume $\operatorname{dim}_{\QQ}V=d+1\geq 2$ with $d\geq 1$, and let $f_0,\cdots,f_d$ be a basis of $V$ over $\QQ$. One can compute such a basis by computing Eisenstein series and multiplying them  \cite{kamal}, \cite{explicitopenimage}*{Algorithm 4.14}. Note that the quotients $f_j/f_i$ are  rational functions of $X_G$. The modular forms $f_0,\cdots,f_d$ define a morphism
\[
\varphi\colon X_G \to \PP_\QQ^d
\]
via $\varphi(P)=[f_0(P),\ldots,f_d(P)]$ for all but finitely many $P$.  Up to an automorphism of $\PP_\QQ^d$, the map $\varphi$ does not depend on the choice of basis. Denote the image of $X_G$ under $\phi$ in $\PP^d_\QQ$ by $C$. Let $I(C)\subseteq \QQ[x_0,\ldots,x_d]$ be the homogeneous ideal of $C$. There is an algorithm to compute a basis for each graded part $I(C)_n$ \cite{explicitopenimage}*{\S 5.3}.

 Let $\calF:=\scrL_k(-E)$ be the invertible sheaf associated to divisor $E$ on $X_G$. The map $\psi_k:M_{k,G} \xrightarrow[]{\sim} \operatorname{H}^0(X_G,\scrL_k )$ restricts to an isomorphism between $V$ and $H^0(X_G,\calF)$ as $\QQ$-vector spaces.  The degree of the invertible sheaf $\scrL_k$ is $k/2\cdot(2g-2)+k/2\cdot r+\lfloor k/4\rfloor v_2+\lfloor k/3 \rfloor \cdot v_3$ where $v_2$ and $v_3$ are the number of elliptic points of $X_G(\CC)$ of order $2$ and $3$. The degree of $\calF$ is given by
\begin{align} \label{E:deg F}
\deg \calF =  \deg \mathscr{L}_k - \sum_{i=1}^r e_i = k/2\cdot (2g-2) + k/2 \cdot r + \lfloor k/4 \rfloor \cdot v_2 +  \lfloor k/3 \rfloor \cdot v_3 - \sum_{i=1}^r e_i.
\end{align}

\subsubsection{Getting the model for $X_G$} \label{L:computing a model}

First, assume $g\geq 3$. In this setting we can choose $k=2$ and the divisor $E=\sum_{i=1}^{r}P_i$, and hence compute the canonical map $\varphi \colon X_G \rightarrow \PP^{g-1}_{\QQ}$ (for more details on the canonical map refer to \cite{actionsoncuspsforms}). If $X_G$ is not geometrically hyperelliptic, then this map is an embedding and $C=\varphi(X_G)$ is a curve isomorphic to $X_G$. The canonical bundle is then normally generated, so the map in \eqref{E:normalgeneration} below is an isomorphism \cite{explicitopenimage}*{\S 5.3.2}.\\
\indent Consider the general case. Note that if $\deg \calF \geq 2g+1$, the Riemann-Roch theorem implies that $\calF$ is very ample, so the map $\varphi$ is an embedding and $C$ is isomorphic to $X_G$. Moreover, $\calF$ is normally generated, i.e., the homomorphism
\begin{equation}\label{E:normalgeneration}
    \eta \colon \QQ[x_0,\ldots, x_d]/I(C) \to \bigoplus_{n\geq 0} H^0(X_G,\calF^{\otimes n})
\end{equation}

defined by $x_i \rightarrow \psi_k(f_i)$ is an isomorphism of $\QQ$-algebras \cite{mumford}. We can choose the even integer $k\geq 2$ large enough and choose the divisor $E=\sum_{i=1}^{r}e_iP_i$ suitably to ensure that $\deg \calF \geq 2g+1$. In this case, the map $\varphi$ is an embedding and $C:=\varphi(X_G)$ is isomorphic to the modular curve $X_G$. 

\subsubsection{Actions on the embedding $C$}\label{SSS:embedding} Let $\calG\subseteq \GL_2(\Zhat)$ be a group such that $G$ is a normal subgroup of $\calG$. One can choose the divisor $E$ suitably so that $\calG$ acts on $V$ and $\operatorname{deg}\scrL_k(-E) \geq 2g+1$ at the same time. The action of $\calG/G$ on $M_{k,G}$ then induces automorphisms of the embedding $C\subseteq \PP_\QQ^d$, defined over $\QQ^{\operatorname{ab}}$, which extend to automorphisms of $\PP_\QQ^d$.

\section{Agreeable Subgroups} \label{S:agreeable}
 In this section, we describe the agreeable subgroups of $\GL_2(\Zhat)$. They were first introduced in \cite{explicitopenimage} and were studied more generally in \cite{openimagenumberfields}. We will mostly follow their exposition.
 \begin{defn}
     We say that a subgroup $\calG$ of $\GL_2(\Zhat)$ is \defi{agreeable} if it is open in $\GL_2(\Zhat)$, satisfies $\det(\calG)=\Zhat^\times$, contains the scalar matrices $\Zhat^{\times} \cdot I$, and the levels of $\calG$ in $\GL_2(\Zhat)$ and $\calG\cap \SL_2(\Zhat)$ in $\SL_2(\Zhat)$ have the same odd prime divisors.
 \end{defn}

Fix an open subgroup $G$ of $\Ggl{\Zhat}$ such that $\operatorname{det}(G)=\Zhat^{\times}$ and $-I\in G$. In general, $G$ will not be an agreeable subgroup. Associated to $G$, there is a unique minimal agreeable subgroup that contains $G$.

\begin{prop}\label{P:agreeable closure}
Let $G$ be an open subgroup of $\GL_2(\Zhat)$ with $\det(G)=\Zhat^\times$.  Let $N$ be the product of primes that divide the level of $[G,G]$ in $\SL_2(\Zhat)$. Define the subgroup
\begin{align}
\calG:= (\ZZ_N^\times \cdot G_N) \times \prod_{\ell\nmid N} \GL_2(\ZZ_\ell)
\end{align}
of $\GL_2(\Zhat)$. Then $\calG$ is the unique minimal agreeable subgroup of $\GL_2(\Zhat)$, with respect to inclusion, that satisfies $G \subseteq \calG$. We have $[G,G]=[\calG,\calG]$ and hence $G$ is a normal subgroup of $\calG$ with $\calG/G$ finite and abelian.
\end{prop}
\begin{proof}
This is Proposition 8.1 in \cite{explicitopenimage}.
\end{proof}
\begin{defn}
    With notation as in Proposition \ref{P:agreeable closure}, we define  the \defi{agreeable closure} of $G$ as the agreeable subgroup $\calG$.
\end{defn}

\begin{remark}
Note that the integer $N$ in Proposition \ref{P:agreeable closure} is even, because the commutator subgroup $[G,G]$ always has even level. The group $\calG$ contains the group $G$ and the scalar matrices $\Zhat^{\times} \cdot I$. The scalar matrices are contained in the center of $\operatorname{GL}_2(\Zhat)$ so $G_N$ and $\Zhat_N^{\times} \cdot G_N$ have the same commutator subgroups. 
\end{remark}

%
%
%

%
%
%
%

\subsection{Constructing the agreeable closure} \label{SS:key construction proof}
The statement of Proposition \ref{P:agreeable closure} suggests that the level of $[G,G]$ needs to be known to compute the agreeable closure $\calG$ of $G$. Unfortunately, computing the commutator subgroup of $G$ can be unfeasible, especially if the level of the group is large or contains large prime factors. However, we can relate the levels of $[\calG,\calG]$ and $\calG \cap \operatorname{SL}_2(\Zhat)$ in $\SL_2(\Zhat)$.

\begin{lemma} \label{L:finally}
\begin{romanenum}
\item\label{L:finally a}
For an odd prime $\ell$, we have $G_\ell=\GL_2(\ZZ_\ell)$ if and only if $\calG_\ell=\GL_2(\ZZ_\ell)$.   
\item \label{L:finally c}
The levels of $[\calG,\calG]$ and $\calG \cap \SL_2(\Zhat)$ in $\SL_2(\Zhat)$ and the level of $\calG$ in $\GL_2(\Zhat)$ have the same odd prime divisors as $N$. 
\end{romanenum}
\end{lemma}
\begin{proof}
This is proven in \cite{explicitopenimage}*{Lemma 8.3}.
\end{proof}

Using the Lemma \ref{L:finally}, one can find the prime divisors of $[G,G]$ from the level of $G \cap \operatorname{SL}_2(\Zhat)$ in $\SL_2(\Zhat)$. The latter is significantly easier to compute.
\begin{lemma}\label{L:oddprime}
Let $G\subseteq \GL_2(\Zhat)$ be an open subgroup with full determinant. Let $B:=G \cap \SL_2(\Zhat)$. Then the level of $[G,G]$ and level of $B$ in $\SL_2(\Zhat)$ have the same odd prime divisors.
\end{lemma}
\begin{proof}
The inclusion $[G,G]\subseteq B$ implies that the level of $B$ divides the level of $[G,G]$.
Let $\calG$ be the agreeable closure of $G$. Since $G\subseteq \calG$, the level of $\calG \cap \Ssl{\Zhat}$ divides the level of $B$. The Lemma  \ref{L:finally} implies that the level of $[\calG,\calG]=[G,G]$ and $\calG \cap \Ssl{\Zhat}$ in $\Ssl{\Zhat}$ have the same odd prime divisors. Therefore, any odd prime $\ell$ that divides the level of $[G,G]=[\calG,\calG]$ also divides the level of $\SL_2(\Zhat) \cap \calG $ and consequently the level of $B$. 
\end{proof}

\begin{prop}
Let $G$ be an open subgroup of $\GL_2(\Zhat)$. Let $\calG$ be the agreeable closure of $G$. Then 
\[\calG=\calG_N \times \prod_{\ell\nmid N} \GL_2(\ZZ_\ell)=(\ZZ_N^\times \cdot G_N) \times \prod_{\ell\nmid N} \GL_2(\ZZ_\ell)\]
where $N$ is the least common multiple of $2$ and the radical of the level of $B=G\cap \operatorname{SL}_2(\Zhat)$ in $\operatorname{SL}_2(\Zhat)$. 
If $G$ has odd level, then so has $\calG$.
\end{prop}
\begin{proof}
By Lemma \ref{L:oddprime}, the levels of $B$ and $[G,G]$ have the same odd prime divisors in $\SL_2(\Zhat)$. From the construction of the agreeable subgroup $ \calG$ the assertion follows.

Since $G \subset \calG$, the level of $\calG$ in $\GL_2(\Zhat)$ divides the level of $G$ in $\GL_2(\Zhat)$. Therefore, if the level of $G$ is an odd integer, so is the level of $\calG$.
\end{proof}

\section{Families of Modular Curves}\label{S:families}

 In this section, we define the families of groups that we use in our classification and collect some results about them. Later in \S\ref{S:twist}, we will show that these families of groups correspond to families of modular curves in a natural way.

Let $\calG$ be an agreeable subgroup of $\GL_2(\Zhat)$. Fix a subgroup $B$ of $\calG$ satisfying $[\calG,\calG] \subseteq B \subseteq \calG \cap \SL_2(\Zhat)$.

\begin{defn}
The \defi{family} of groups associated to the pair $(\calG,B)$ is the set $\scrF(\calG,B)$ of open subgroups $H$ of $\calG$ that satisfy $\det(H)=\Zhat^\times$ and $H\cap \SL_2(\Zhat)=B$.
\end{defn}

Assume that $\scrF(\calG,B)$ is nonempty. Fix a group $G\in \scrF(\calG,B)$. We know that $G$ is an open subgroup of $\calG$. Since $[\calG,\calG] \subset G$, we have that $G$ is a normal subgroup of $\calG$, and the quotient $\calG/G$ is finite and abelian. Consider any continuous homomorphism $\gamma\colon \Zhat^\times \to \calG/G$. Define the group
\[
G_\gamma:=\{g \in \calG :  g\cdot G = \gamma(\det g) \}.
\]

\begin{lemma}\label{L:twistgroup} \cite{explicitopenimage}*{Lemma 14.2} \label{familygroups}
With notation as above, the set $\scrF(\calG,B)$ consists of the groups $G_\gamma$ with $\gamma\colon \Zhat^\times \to \calG/G$ a continuous homomorphism.
\end{lemma}
\begin{proof}
\indent First take any $\gamma$.  We have $G_\gamma \cap \SL_2(\Zhat) = G\cap \SL_2(\Zhat)=B$.    The natural map $(\calG \cap \SL_2(\Zhat))/B \to \calG/G$ is an isomorphism since $G\cap \SL_2(\Zhat) = B$ and $\det(G)=\Zhat^\times$.  Using this isomorphism, we find that $\det(G_\gamma)=\Zhat^\times$.   Therefore, $G_\gamma \in \scrF(\calG,B)$.

\indent Conversely, take any $H \in \scrF(\calG,B)$.  The quotient map $H\to \calG/G$ induces a homomorphism $f\colon H/(H\cap \SL_2(\Zhat)) \to \calG/G$ since $H\cap \SL_2(\Zhat)=B=G\cap \SL_2(\Zhat)$.   Let $\gamma\colon \Zhat^\times \to \calG/G$ be the homomorphism obtained by composing the inverse of the determinant map $H/(H\cap \SL_2(\Zhat)) \xrightarrow{\sim} \Zhat^\times$ with $f$.  For each $h\in H$, we have $h \cdot G = \gamma(\det h)$.  Therefore, $H \subseteq G_\gamma$.  Since $H$ and $G_\gamma$ both have full determinant and have the same intersection with $\SL_2(\Zhat)$, we conclude that $H=G_\gamma$.
\end{proof}

 


Let $N$ be the least common multiple of the levels of $\calG$ in $\GL_2(\Zhat)$ and $B$ in $\SL_2(\Zhat)$. We define $N_1\colonequals N$ if $N$ is odd and $N_1\colonequals \operatorname{lcm}(N,8)$ if $N$ is even. Then

\begin{thm} \label{T:representativeinfam}
Let $S$ be the $2$-power torsion subgroup of $\ZZ_N^{\times}$. Then the following are equivalent:
\begin{enumerate}
\item There is an open subgroup $G\subseteq \calG$ with $G\cap \SL_2(\Zhat) = B$ and $\operatorname{det}(G)=\Zhat^{\times}$.
\item There is a homomorphism $\beta\colon S \rightarrow \calG_N/B_N$ such that $\operatorname{det}(\beta(a))=a$ for all $a\in S$.
\item There is a homomorphism $\beta\colon S \rightarrow \calG(N_1)/B(N_1)$ such that $\operatorname{det}(\beta(a))\equiv a \; (\operatorname{mod }N_1)$ for all $a\in S$.
\end{enumerate}
Moreover, if a group $G$, as in $(1)$ exists, then there is such a group whose level divides a power of $2$ times $N$.
\end{thm}
\begin{proof}
    This is a special case of \cite{openimagenumberfields}*{Theorem 4.5} with $U=\Zhat^{\times}$.
\end{proof}

    \begin{thm}\label{T:indexandgenus}
        Take $G,H \in \scrF(\calG,B)$. Then
        \begin{enumerate}
            \item $X_G$ and $X_H$ have the same genus.
            \item $[\GL_2(\Zhat): G] =[\GL_2(\Zhat): H]$.
        \end{enumerate}
    \end{thm}
    \begin{proof}
(1) First note that $\Gamma_G =\Gamma_H$. Since $X_H(\CC)$ and $X_G(\CC)$ are both isomorphic to $\calX_{\Gamma_G}$ as Riemann surfaces. The assertion follows.\\
(2) Since $G$ and $H$ both have full determinant, we have that $[\GL_2(\Zhat): G]=[\SL_2(\Zhat): B]=[\GL_2(\Zhat): H]$. Therefore, $G$ and $H$ have the same index in $\GL_2(\Zhat)$.

            
    \end{proof}

For an open subgroup $G\subseteq \GL_2(\Zhat)$, there is a natural choice of family for $G$.

\begin{cor}\label{C:aggcontainsgroup}
Let $G$ be an open subgroup of $\GL_2(\Zhat)$ with full determinant, and let $B= G \cap \operatorname{SL}_2(\Zhat)$.  Let $\calG$ be the agreeable closure of $G$. Then $G\in \scrF(\calG,B)$.
\end{cor}
\begin{proof}
We have $G \subseteq \calG$. Since the commutator subgroups of $G$ and $\calG$ agree (Proposition \ref{P:agreeable closure}), we have $[\calG,\calG]=[G,G] \subseteq B \subseteq \calG \cap \SL_2(\Zhat)$, implying that $G\in \scrF(\calG,B)$.
\end{proof}
\begin{remark}
An open subgroup $G$ may lie in more than one family. However, Corollary \ref{C:aggcontainsgroup} suggests a canonical choice of family that contains $G$, i.e., the family $\scrF(\calG,B)$ where $\calG$ is the agreeable closure of $G$ and $B=G\cap \Ssl{\Zhat}$. We have a description of $\calG$ that depends on $G$ and $B$. In particular, one can easily compute the group $\calG$ and subsequently identify the family $\scrF(\calG,B)$ which contains $G$.
\end{remark}
\subsection{Translating continuous homomorphisms:}\label{SS:translating}
Let $\scrF(\calG,B)$ be a non-empty family of modular curves. Fix a group $G\in\scrF(\calG,B)$. For any continuous homomorphism $\gamma:\Zhat^{\times} \rightarrow \calG/G$, there is a twisted group $G_\gamma \in \scrF(\calG,B)$ defined as in Lemma \ref{L:twistgroup}. Conversely, if $H\in \scrF(\calG,B)$, Lemma \ref{L:twistgroup} shows that $H=G_{\gamma_H}$ for a unique continuous homomorphism $\gamma_H:\Ztimes\rightarrow \calG/G$. 

Let $G_\gamma\in \scrF(\calG,B)$ where $\gamma: \Ztimes \rightarrow \calG/G$ is a continuous homomorphism and let $W=\calG \cap \SL_2(\Zhat)  $. The natural map $W/B\rightarrow \calG/G_\gamma$ is an isomorphism. Let $f_\gamma: \calG/G_\gamma \xrightarrow{\sim} W/B$ denote the inverse of this isomorphism. More specifically, if $gG_\gamma \in \calG/G_\gamma$ with $\det(g)=d$, then $f_\gamma(gG_\gamma)=tgB$ where $t\in G_\gamma$ is a matrix that satisfies $\det(t)=d^{-1}$. The isomorphism $f_\gamma$ is independent of the choices of $t$. When $\gamma$ is the trivial homomorphism, we define $f:=f_\gamma$.


Let $\alpha: \Ztimes \rightarrow \calG/G_\gamma$ be a continuous homomorphism. Twisting $G_\gamma$ with $\alpha$ one obtains a group $(G_\gamma)_\alpha$. Since $(G_\gamma)_\alpha$ lies in $\scrF(\calG,B)$, it is of the form $G_\psi$ for some continuous homomorphism $\psi: \Ztimes \rightarrow \calG/G$. We characterize $\psi$ in terms of $\gamma$ and $\alpha$. 
\begin{lemma}\label{L:translate}
    With notation as above, we have $(G_\gamma)_\alpha=G_\psi$ where $\psi:\Zhat^{\times}\rightarrow \calG/G$ is defined by $\psi(d):= \gamma(d)f^{-1}(f_\gamma(\alpha(d)))$ for $d\in \Ztimes$.
\end{lemma}
\begin{proof}
    By definition, we have
\[
G_\psi:=\{g \in \calG :  g\cdot G = \psi(\det g) \}
\]
and 
\[
(G_\gamma)_\alpha:=\{g \in \calG :  g\cdot G_\gamma = \alpha(\det g) \}.
\]
We check that these two groups are the same. Let $g\in (G_\gamma)_\alpha$ and let $d=\det(g)$. We have
\begin{align*}
    \psi(d)=\gamma(d)\cdot f^{-1}(f_\gamma(\alpha(d)))=\gamma(d)\cdot (tg)G
\end{align*}
where $t\in G_\gamma$ with $\det(t)=d^{-1}$. Note that, $t$ satisfies $\gamma(\det(t))=t G$ so we have $\gamma(d)=\gamma(\det(t^{-1})) =t^{-1}G$. Therefore,
\[\psi(d)=t^{-1}G\cdot tgG=gG\]
showing that $g\in G_\psi$ and $(G_\gamma)_\alpha\subseteq G_\psi$. Since $(G_\gamma)_\alpha \cap \SL_2(\Zhat)=G_\psi\cap \SL_2(\Zhat)$ and both groups have full determinant, we have $(G_\gamma)_\alpha=G_\psi$.

\end{proof}
By Lemma \ref{L:translate}, we have a map $\tau_\gamma:\operatorname{Hom}_{\operatorname{cts}}(\Ztimes,\calG/G_\gamma) \rightarrow \operatorname{Hom}_{\operatorname{cts}}(\Ztimes,\calG/G)$ defined by $\tau_\gamma(\alpha)(d)=\gamma(d) \cdot f^{-1}(f_\gamma(\alpha(d)))$ and satisfying $G_{\tau_\gamma(\alpha)}=(G_\gamma)_\alpha$, where $\operatorname{Hom}_{\operatorname{cts}}(\Ztimes,\calG/G)$ denotes the continuous homomorphisms from $\Ztimes$ to $\calG/G$.

Let $\gamma_1$ and $\gamma_2$ be two continuous homomorphisms $\Ztimes \rightarrow \calG/G$. Let $\gamma_1\gamma_2$ be the product homomorphism. The Lemma \ref{L:translate} implies that 
\[G_{\gamma_1\gamma_2}=(G_{\gamma_1})_{\tau_{\gamma_1}^{-1}(\gamma_1\gamma_2)}. \]
Moreover,
\[
\tau_{\gamma_1}^{-1}(\gamma_1\gamma_2)=f_{\gamma_1}^{-1}\circ f\circ\gamma_2.
\]
Thus multiplication by $\gamma_1^{-1}$ is followed by transport from $\calG/G$ to $\calG/G_{\gamma_1}$ through their identifications with $W/B$. If we let $\calG/G \simeq \ZZ/p_1^{r_1}\ZZ \times \ldots \times \ZZ/p_s^{r_s}\ZZ$ for some primes $p_1,\ldots,p_s$, we see that any continuous homomorphism $\delta: \Ztimes \rightarrow \calG/G$ can be written as a product of homomorphisms whose image is a cyclic group of prime power order. Using this fact, for certain computations we will restrict our attention to continuous homomorphisms $\Ztimes\rightarrow \calG/G$ that factor through a cyclic factor $\ZZ/p^r\ZZ$ of $\calG/G$, for a prime number $p$ and a positive integer $r$.

\section{Finiteness of Agreeable Subgroups} \label{S:finiteness}

Fix a non-negative integer $g$. In this section, we recall a result from \cite{openimagenumberfields} that there are finitely many agreeable subgroups up to conjugacy, of genus less than or equal to $g$. We will also describe a method for computing all such agreeable subgroups. We first start by making some observations.

\indent Let $\calG$ be an agreeable subgroup of genus at most $g$. Let $H:=\calG \cap \operatorname{SL}_2(\Zhat)$; it is an open subgroup in $\operatorname{SL}_2(\Zhat)$. Let $N$ be the level of $H$ in $\SL_2(\Zhat)$. The associated congruence subgroup $\Gamma_G:=H \cap \SL_2(\ZZ)$ agrees with the congruence subgroup of level $N$ consisting of elements in $\operatorname{SL}_2(\ZZ)$ whose image modulo $N$ lies in $H$ modulo $N$. Similarly we have that $-I \in \Gamma_G$, and $\Gamma_G$ has genus at most $g$. Applying the finiteness theorem of Cox and Parry \cite{cox-parry-genus}*{Corollary 4.7} to $\PSL_2(\ZZ)$ and taking inverse images in $\SL_2(\ZZ)$, we see that there are only finitely many congruence subgroups of $\operatorname{SL}_2(\ZZ)$ of genus at most $g$ containing $-I$ (see also Thompson \cite{thompson-finiteness} for finiteness up to conjugacy in $\PSL_2(\RR)$). The database of Cummins and Pauli \cite{cpdata} lists all such subgroups of genus at most $24$ up to conjugacy in $\SL_2(\ZZ)$, with the counts given in \cite{MR2016709}*{Theorem 2.8}.\\
\indent In the proof, we will reverse this process and explain how to obtain the finitely many agreeable subgroups up to genus $g$ arising from a congruence subgroup $\Gamma$ of genus at most $g$. We start with a lemma.

\begin{thm}\label{T:finmanyagg}
There are finitely many agreeable subgroups of $\operatorname{GL}_2(\Zhat)$ up to conjugacy with genus at most $g$. 
\end{thm} 
\begin{proof}
Fix a genus $g$ and fix a congruence subgroup $\Gamma$ that has genus at most $g$ and contains $-I$. There are finitely many such congruence subgroups. Consider the corresponding subgroup of $\operatorname{SL}_2(\Zhat)$ which we call $H$. The level of $H$ is equal to the level of $\Gamma$, which we call $N$. In particular $H$ is the subgroup of $\operatorname{SL}_2(\Zhat)$ whose image modulo $N$ is equal to $\Gamma$ modulo $N$. \\
\indent We first determine how the level of an agreeable subgroup $G\subseteq \operatorname{GL}_2(\Zhat)$ such that $G\cap \operatorname{SL}_2(\Zhat)= H$ relates to the level of $H$.  Let $N_1:=2\cdot \operatorname{lcm}(N,12)$.
\begin{lemma}\label{L:agglevel}
Any agreeable subgroup $G\subseteq \Ggl{\Zhat}$ with $G \cap \Ssl{\Zhat} = H$ has level dividing $N_1$. 
\end{lemma}
\begin{proof}
See \cite{openimagenumberfields}*{Lemma 5.1}.

\end{proof}
Let $\Gamma$ and $H$ be as above. Assume $G\subseteq \operatorname{GL}_2(\Zhat)$ is an agreeable subgroup with $G \cap \operatorname{SL}_2(\Zhat) = H$. By Lemma \ref{L:agglevel}, the group $G$ has level dividing $N_1$, so $G$ corresponds to a subgroup $\bar{G}$ of $\operatorname{GL}_2(\ZZ/N_1\ZZ)$ such that $\bar{G}\cap \operatorname{SL}_2(\ZZ/N_1\ZZ) = \bar{H}$ where $\bar{H}$ denotes the reduction of $H$ to modulo $N_1$. There are only finitely many such subgroups of $\operatorname{GL}_2(\ZZ/N_1\ZZ)$, so only finitely many agreeable subgroups of $\operatorname{GL}_2(\Zhat)$ arise from a fixed congruence subgroup $\Gamma$. Since there are finitely many congruence subgroups $\Gamma$ of genus at most $g$ that contain $-I$, we conclude that there are finitely many agreeable subgroups with genus at most $g$.
\end{proof}

\subsection{Computing agreeable subgroups}\label{SS:computeagg}

We explain how to explicitly compute all agreeable subgroups of $\GL_2(\Zhat)$ of genus at most $g$. The effective level bound of Cox and Parry \cite{cox-parry-genus}*{(4.22)} (see also Cummins and Pauli \cite{MR2016709}*{Proposition 2.1}) allows us first to enumerate all congruence subgroups $\Gamma$ containing $-I$ and of genus at most $g$, by searching the finite groups $\SL_2(\ZZ/N\ZZ)$ for the bounded set of levels $N$. For each such $\Gamma$, let $N_1$ be as in the proof of \ref{T:finmanyagg}. We directly search in $\operatorname{GL}_2(\ZZ/N_1\ZZ)$ for subgroups $\bar{G}$ with $(\ZZ/N_1\ZZ)^{\times}\cdot I\subseteq \bar{G}$, $\operatorname{det}(\bar{G})=(\ZZ/N_1\ZZ)^{\times}$, $-I \in \bar{G}$ and $\bar{G}\cap \operatorname{SL}_2(\ZZ/N_1\ZZ) = \bar{H}$. These groups give rise to finitely many, potentially agreeable, subgroups $G\subseteq \operatorname{GL}_2(\Zhat)$ such that $G \cap \operatorname{SL}_2(\Zhat) = H$. For each such $G$, we then check if it is an agreeable subgroup, i.e., if its level in $\operatorname{GL}_2(\Zhat)$ has the same odd prime divisors as the level of $G \cap \operatorname{SL}_2(\Zhat)$ in $\SL_2(\Zhat)$.


We are ready to state the first part of Theorem \ref{T:main families1} in terms of families of groups.

\begin{thm} \label{T:main familiesgroups}
Fix a non-negative integer $g$. There are only finitely many families of groups of genus $g$.
\end{thm}
\begin{proof}
By Theorem \ref{T:finmanyagg}, there are only finitely many agreeable subgroups of $\Ggl{\Zhat}$ up to conjugacy of genus at most $g$. We denote a set of representatives by $\mathscr{A}_g$. For each agreeable subgroup $\calG \in \mathscr{A}_g$, the groups $[\calG,\calG]$ and $\calG \cap \operatorname{SL}_2(\Zhat)$ are open subgroups of $\operatorname{SL}_2(\Zhat)$, hence there are only finitely many subgroups $B$ of $\calG$ such that $[\calG,\calG] \subseteq B \subseteq \calG \cap \Ssl{\Zhat}$. Combining the finiteness of agreeable subgroups $\calG$ of genus at most $g$ and the finitely many subgroups $B$ arising from $\calG$, we obtain a finite set of pairs $(\calG_i,B_i)$ and associated families $\scrF(\calG_i,B_i)$. \\
\indent Let $G$ be an open subgroup of $\GL_2(\Zhat)$ of genus $g$ with full determinant and $-I \in G$. By Corollary \ref{C:aggcontainsgroup}, $G$ lies in the family $\scrF(\calG,B)$ where $B=G \cap \Ssl{\Zhat}$ and $\calG$ is the agreeable closure of $G$. The agreeable closure $\calG$ has genus at most $g$ and so, after simultaneous conjugation of $G$, $\calG$ and $B$, it is in the finite set $\mathscr{A}_g$. Hence, $\scrF(\calG,B)$ is among the finite set of families $\scrF(\calG_i,B_i)$ and these families cover the set of all open subgroups $G \subseteq \GL_2(\Zhat)$ of genus $g$ with full determinant and $-I\in G$.
\end{proof}

\section{Twisting Modular Curves} \label{S:twist}

We have stated in \S \ref{SS:famsfirst} that a family of groups $\scrF(\calG,B)$ corresponds to a family of twists of modular curves. In this section, we will describe the spaces of modular forms $M_{k,G}$ as we vary $G$ in a family $\scrF(\calG,B)$ and, consequently, describe how to twist the modular curves $X_G$. 
\subsection{Twisting spaces of modular forms:} Fix a nonempty family $\scrF(\calG,B)$ and a group $G\in \scrF(\calG,B)$. Let $X_G$ be the modular curve associated to $G$, and let $\pi_{\calG/G}: X_G \rightarrow X_{\calG}$ be the morphism coming from the inclusion $G \subseteq \calG$. We start with a definition:
\begin{defn}
A \defi{$\calG-$twist} of $(X_G,\pi_{\calG/G})$ is a pair $(Y,\pi)$ where $Y$ is a curve over $\QQ$, with a morphism $\pi:Y\rightarrow X_{\calG}$ defined over $\QQ$, such that there is an isomorphism $f:(X_G)_{\QQ^{\operatorname{ab}}} \rightarrow (Y)_{\QQ^{\operatorname{ab}}}$ that satisfies $\pi \circ f = \pi_{\calG/G}$. We call the pairs $(X_G,\pi_{\calG/G})$ and $(Y,\pi)$ \defi{isomorphic}, if there is such an $f$ defined over $\QQ$. 
\end{defn}

The group $G$ is a normal subgroup of $\calG$. The latter acts on $M_{k,G}$ for all $k> 0$ via the $\ast$ action of Lemma \ref{L:star action}. By definition, $G$ acts trivially on $M_{k,G}$ so there is an action of $\calG/G$ on $M_{k,G}$. Consequently, $\calG/G$ acts on $X_G$ by Definition \ref{D:modcurve}. The degree of $\pi_{\calG/G}$ is $|\calG/G|$, therefore we have $\operatorname{Aut}(X_G/X_\calG)= \calG/G$ where $\operatorname{Aut}(X_G/X_\calG)$ is the group of automorphisms $f$ of the curve $X_G$ that satisfy $\pi_{\calG/G} \circ f = \pi_{\calG/G}$. Note that these are \textit{modular automorphisms}, and since there is a natural isomorphism $\calG/G \simeq (\calG \cap \SL_2(\Zhat))/B$, the automorphisms in $\operatorname{Aut}(X_G/X_\calG)$ are defined over $\QQ$. 
In this identification, $gG$ corresponds to the inverse of the automorphism induced by its $\ast$-action on the graded ring of modular forms. This is a group isomorphism since $\calG/G$ is abelian.

Let $\gamma : \Zhat^{\times} \rightarrow \calG/G$ be a continuous homomorphism. By precomposing with the cyclotomic character $\chi_{\operatorname{cyc}}$ we obtain a homomorphism 
\[\xi := \gamma \circ \chi_{\texttt{cyc}}: \operatorname{Gal}(\QQ^{\operatorname{ab}}/\QQ) \rightarrow \calG/G \cong  \operatorname{Aut}(X_G/X_\calG) .
\]

\noindent In particular, $\xi$ is a 1-cocycle of $X_G$. 

\begin{lemma}
There is a bijection between isomorphism classes of $\calG-$twists of $X_G$ and
\[
\operatorname{H}^1(\Gal(\QQ^{\operatorname{ab}}/\QQ),\operatorname{Aut}(X_G/X_\calG)).
\]
\end{lemma}
\begin{proof}
Let $(X,\pi_X):=(X_G,\pi_{\calG/G})$, and let $(Y,\pi_Y)$ be a $\calG-$twist of $(X,\pi_X)$. So there is an isomorphism $f:X_{\QQ^{\operatorname{ab}}}\rightarrow Y_{\QQ^{\operatorname{ab}}}$ such that $\pi_Y \circ f =\pi_X$. Let $\xi: \operatorname{Gal}(\QQ^{\operatorname{ab}}/\QQ)\rightarrow \operatorname{Aut}(X_G/X_\calG)$ be defined as $\xi_\sigma=f^{-1}\circ \sigma(f)$. One can check that $\xi$ is a 1-cocycle. We have $\pi_X \circ \xi_\sigma=\pi_X \circ f^{-1} \circ \sigma(f)=\pi_Y \circ \sigma(f)=\sigma(\pi_Y \circ f)=\sigma(\pi_X)=\pi_X$, because $\pi_X$ and $\pi_Y$ are defined over $\QQ$. We define the map $\lambda$ such that $\lambda((Y,\pi_Y))=[\xi]$.

We first show that $\lambda$ is well defined and does not depend on the choice of the isomorphism. Let $g:X_{\QQ^{\operatorname{ab}}} \rightarrow
 Y_{\QQ^{\operatorname{ab}}}$ be another such isomorphism. Then $ f^{-1} \circ \sigma(f)$ and $ g^{-1} \circ \sigma(g)$ are cohomologous, so the class of the cocycle is well defined.
 
 Let $Y$ and $Z$ be two curves that are isomorphic over $\QQ$ with a map $h:Y \rightarrow Z$ such that $\pi_Z \circ h = \pi_Y$. Then we have an isomorphism $h\circ f : X_{\QQ^{\operatorname{ab}}} \rightarrow Z_{\QQ^{\operatorname{ab}}}$ satisfying $\pi_Z \circ h \circ f =\pi_X$. Looking at associated cocycles, we have $(h\circ f)^{-1} \circ \sigma(h\circ f)=f^{-1}\circ h^{-1}\circ \sigma(h)\circ \sigma(f)=f^{-1}\circ \sigma(f)$ because $h$ is defined over $\QQ$. So, the $[\xi]$ is independent of the $\QQ$ isomorphism class of $Y$.
 
 Let $\xi_1$ and $\xi_2$ be two cocycles that are cohomologous corresponding to $Y_1$ and $Y_2$. This means that there is $T\in \operatorname{Aut}(X_G/X_\calG)$ such that $\xi_1(\sigma)=T^{-1}\circ \xi_2(\sigma)\circ \sigma(T)$. Then $f_2 \circ T \circ f_1^{-1}: Y_1 \rightarrow Y_2$ is an isomorphism defined over $\QQ$, proving the injectivity of $\lambda$.

 To show surjectivity, let $[\xi]$ be a class in $H^1(\Gal(\QQ^{\operatorname{ab}}/\QQ),\operatorname{Aut}(X_G/X_\calG))$. By Galois descent, we have a smooth, projective, geometrically integral curve $Y$ over $\QQ$ with an isomorphism $f:X_{\QQ^{\operatorname{ab}}}\rightarrow Y_{\QQ^{\operatorname{ab}}}$ where $\sigma \rightarrow f^{-1}\circ \sigma(f)$ is cohomologous to $\xi$ \cite[Chapter III, Proposition 5]{Serre}. Precisely this means that there is $T\in \operatorname{Aut}(X_G/X_\calG)$ such that  $f^{-1}\circ \sigma(f)=T^{-1}\circ \xi_\sigma \circ \sigma(T)$. Set $g=f\circ T^{-1}$. We have $g^{-1}\circ \sigma(g)=\xi_\sigma$. Define $\pi_Y:=\pi_X \circ g^{-1}$. Note that we have $\pi_X \circ \xi_\sigma=\pi_X$. Let $P \in Y_{\QQ^{\operatorname{ab}}}$ let $Q=g^{-1}(\sigma^{-1}(P))$ and $P':=\sigma^{-1}(P)$. We have
 \[
 \sigma(\pi_Y)(P)=\sigma(\pi_X(g^{-1}(P')))
 \]
$\pi_X$ is defined over $\QQ$ so
\[
\sigma(\pi_X(g^{-1}(P')))=\pi_X(\sigma(g^{-1}(P')))=\pi_X(\sigma(g^{-1})(P))
\]
which is equal to 
\[
\pi_X(\xi_\sigma^{-1}\circ g^{-1}(P))=\pi_X(g^{-1}(P))=\pi_Y(P).
\]
Hence $\pi_Y=\sigma(\pi_Y)$ and $(Y,\pi_Y)$ is the associated $\calG$-twist of $X$.



\end{proof}

Let $\gamma : \Zhat^{\times} \rightarrow \calG/G$ be a continuous homomorphism and let $\xi: \operatorname{Gal}(\QQ^{\operatorname{ab}}/\QQ) \rightarrow \calG/G =  \operatorname{Aut}(X_G/X_\calG) $  be the associated cocycle by precomposing with the cyclotomic character. Twisting via this cocycle, we get a $\calG$-twist $(\left(X_G\right)_{\xi},(\pi_{\calG/G})_{\xi})$. We prove that this pair is isomorphic to $(X_{G_\gamma},\pi_{\calG/G_\gamma})$, where $G_\gamma$ is the group defined in \S \ref{S:families}.

\indent Observe that $G\cap \Ssl{\Zhat} = G_\gamma \cap \Ssl{\Zhat}=B$. Put $\Gamma_B:=B\cap\SL_2(\ZZ)$ and define $M_{k,B}:=M_k(\Gamma_B,\QQ^{\operatorname{ab}})$, where each form has Fourier coefficients in a finite cyclotomic extension. Lemma \ref{L:tensor agreement} gives the following equalities for even $k$; for odd $k$, all three spaces vanish since $-I\in B$:
\[
M_{k,B}=M_{k,G}\otimes_\QQ \mathbb{Q}^{\operatorname{ab}} =  M_{k,G_\gamma}\otimes_\QQ \mathbb{Q}^{\operatorname{ab}}.
\]
\indent Let $g\in \calG$ be an element and let $\sigma \in \operatorname{Gal}(\mathbb{Q}^{\operatorname{ab}}/\QQ)$ be such that $\chi_{\operatorname{cyc}}(\sigma)=\operatorname{det}(g)$. The group $\calG$ acts on $M_{k,G}\otimes_\QQ \mathbb{Q}^{\operatorname{ab}}$ where $g$ sends $f\otimes c$ to $(f \ast g) \otimes \sigma(c)$. 


\indent Denote by $cf$, the element $f\otimes c \in M_{k,G}\otimes_\QQ \mathbb{Q}^{\operatorname{ab}}$. We define twisted action of $\operatorname{Gal}(\mathbb{Q}^{\operatorname{ab}}/\QQ)$ on $M_{k,G}\otimes_\QQ \mathbb{Q}^{\operatorname{ab}}$ by $\sigma \bullet (cf) :=\sigma(c)(\xi_\sigma(f))$ where $\xi_\sigma(f)$ denotes the action of $\calG/G$ on $M_{k,G}$ via the cocycle, equivalently inverse pullback by the corresponding curve automorphism. For each $k\geq0$, we define the twisted space $(M_{k,G})_\xi$ by
\begin{align*}
(M_{k,G})_\xi = \{f \in M_{k,G}\otimes_\QQ \mathbb{Q}^{\operatorname{ab}} \;\vert \; \sigma \bullet f = f \; \text{for all}\; \sigma \in \operatorname{Gal}(\mathbb{Q}^{\operatorname{ab}}/\QQ)\}.
\end{align*}


Restricting the action of $\calG$ to $G$ and $G_{\gamma}$, we obtain the induced actions of $G/B$ and $G_\gamma/B$ on $ M_{k,G}\otimes_\QQ \mathbb{Q}^{\operatorname{ab}}$. Composing these with the isomorphisms 
\[\varphi_1\colon \operatorname{Gal}(\mathbb{Q}^{\operatorname{ab}}/\QQ)\rightarrow \Zhat^{\times} \xrightarrow{\operatorname{det}^{-1}}  G/B\] 
and
\[
\varphi_2\colon \operatorname{Gal}(\mathbb{Q}^{\operatorname{ab}}/\QQ)\rightarrow \Zhat^{\times} \xrightarrow{\operatorname{det}^{-1}} G_\gamma/B
\]
we get two different Galois actions of $\Gal(\QQ^{\operatorname{ab}}/\QQ)$ on $M_{k,G}\otimes_\QQ \mathbb{Q}^{\operatorname{ab}}$ and hence on $M_{k,B}$. 

\indent It is important that the action $\bullet$ of $\operatorname{Gal}(\mathbb{Q}^{\operatorname{ab}}/\QQ)$ and the action of $G_\gamma/B$ on $M_{k,G}\otimes_\QQ \mathbb{Q}^{\operatorname{ab}}$ are compatible in the following sense. 
\begin{lemma}
If $\sigma \in \operatorname{Gal}(\mathbb{Q}^{\operatorname{ab}}/\QQ)$ and $g\in G_\gamma$ with $\operatorname{det}(g) = \chi_{\operatorname{cyc}}(\sigma)$ then $\sigma \bullet f = f\ast g$ for all $f \in M_{k,G}\otimes_\QQ \mathbb{Q}^{\operatorname{ab}}$. 
\end{lemma}
\begin{proof}
Let $cf:= f\otimes c \in M_{k,G}\otimes_\QQ \mathbb{Q}^{\operatorname{ab}}$. We have $\sigma \bullet (cf)=\sigma(c)(\xi_\sigma(f))=(\xi_\sigma(f)) \otimes \sigma(c)$. Let $g$ be as in the statement of the lemma. Then we have that 
\[\xi_\sigma=\gamma(\operatorname{det}(g))=gG\]
in $\calG/G$. \\
\indent Hence $(cf)\ast g = \sigma(c) (f \ast g)=f \ast g \otimes \sigma(c) =  \xi_\sigma (f) \otimes \sigma (c)=\sigma\bullet(cf)$.
\end{proof}

\begin{thm}\label{T:mkg}
Let $\gamma:\Zhat^{\times}\rightarrow \calG/G$ be a continuous homomorphism. Let $\xi=\gamma\circ \chi_{\operatorname{cyc}}: \Gal(\QQ^{\operatorname{ab}}/\QQ)\rightarrow \calG/G$ be the associated cocycle. Then $(M_{k,G})_\xi=M_{k,G_\gamma}$.
\end{thm}
\begin{proof}
Let $f\in M_{k,G_\gamma}$.  Then for all $\sigma\in \Gal(\QQ^{\operatorname{ab}}/\QQ)$ and $g\in G_\gamma$ such that $\det(g)=\chi_{\operatorname{cyc}}(\sigma)$, we have $\sigma \bullet f = f\ast g = f$, which implies that $f \in (M_{k,G})_\xi$.\\
\indent For the converse, let $f\in (M_{k,G})_\xi$. Then for all $g\in G_\gamma$ and $\sigma\in \Gal(\QQ^{\operatorname{ab}}/\QQ)$ such that $\det(g)=\chi_{\operatorname{cyc}}(\sigma)$, we have $f\ast g = \sigma \bullet f=f$. Hence $f\in M_{k,G_\gamma}$.
\end{proof} 
The action $\bullet$ on $M_{k,G}\otimes_\QQ \QQ^{\operatorname{ab}}$ induces an action of $\operatorname{Gal}(\mathbb{Q}^{\operatorname{ab}}/\QQ)$ on the $\QQ^{\operatorname{ab}}$-algebra $\bigoplus _{k=0}^{\infty}M_{k,G} \otimes_\QQ \QQ^{\operatorname{ab}}$. We define
\[
(R_G)_\xi:= \bigoplus_{k=0}^{\infty}(M_{k,G})_\xi.
\]
We immediately get the following:
\begin{cor}\label{C:twist} With notation as above, the pair $(X_{G_\gamma},\pi_{\calG/G_\gamma})$ is isomorphic to $((X_{G})_\xi,(\pi_{\calG/G})_\xi)$.
\end{cor}
\begin{proof}
    By definition, we have $X_{G_\gamma}=\operatorname{Proj}(\bigoplus_{k=0}^{\infty}M_{k,G_\gamma})$, where the nonconstant map $\pi_{\calG/G_\gamma}$ is induced by the inclusions $M_{k,\calG}\subseteq M_{k,G_\gamma}$. The $\calG$-twist $((X_G)_\xi,(\pi_{\calG/G})_\xi)$ is isomorphic to $(\operatorname{Proj}((R_G)_\xi),\pi)$, where $\pi$ is induced by the inclusions $M_{k,\calG}\subseteq (M_{k,G})_\xi$. Theorem \ref{T:mkg} implies that these two pairs are isomorphic. 
\end{proof}

We conclude from Corollary \ref{C:twist} that our family of groups $ \scrF(\calG,B)$ in fact corresponds to a family of abelian twists of modular curves, i.e. it consists of curves of the form $X_{G_\gamma}\simeq(X_G)_\xi$.

We denote by $\scrF_g$ the finite collection, up to conjugacy, of nonempty families $\scrF(\calG,B)$ with $\calG\in\scrA_g$, $-I\in B$, and $B\cap\SL_2(\ZZ)$ of genus at most $g$. This collection has been explicitly determined for $g=24$ \cite{repo}.

\subsubsection{Getting a basis for $M_{k,G_{\gamma}}$} \label{SS:getting a basis}Assume that, using the methods described in \S\ref{S:Models}, we have an explicit basis $\mathcal{B}:=\{ f_0,\cdots, f_d\}$ for $M_{k,G}$. In particular, one can use Eisenstein series of weight 1 to compute such a basis, cf. \cite{kamal} and \cite{explicitopenimage}*{Algorithm 4.14}. The group $\calG/G$ is finite and abelian. It also acts on $M_{k,G}$ via $\QQ$-linear automorphisms. We can compute the action of any $gG\in \calG/G$ on the basis $\mathcal{B}$ and get a matrix in $\operatorname{GL}_{d+1}(\mathbb{Q})$. Hence, for any 1-cocycle $\xi:\operatorname{Gal}(\mathbb{Q}^{\operatorname{ab}}/\QQ) \rightarrow \calG/G$, we get a related cocycle 
\[
\bar{\xi}:\operatorname{Gal}(\mathbb{Q}^{\operatorname{ab}}/\QQ) \rightarrow \operatorname{GL}_{d+1}(\mathbb{Q}).
\]
Write $\mathbf f=(f_0,\ldots,f_d)^{\mathsf T}$. The coefficient row of a form $h=\sum_{i=0}^d w_i f_i$ with $w_i\in\QQ^{\operatorname{ab}}$ is $w=(w_0,\ldots,w_d)$, so $h=w\mathbf f$. For $g_\sigma G=\xi_\sigma$, our convention is
\[
f_i\ast g_\sigma=\sum_{j=0}^d\bar\xi(\sigma)_{ij}f_j.
\]
Thus the twisted action on a coefficient row $w$ is $\sigma\bullet w=\sigma(w)\bar\xi(\sigma)$.

By Hilbert 90, $\operatorname{H}^1(\operatorname{Gal}(\mathbb{Q}^{\operatorname{ab}}/\QQ),\operatorname{GL}_{d+1}(\mathbb{Q}^{\operatorname{ab}}))$ is trivial so there exists a matrix $A\in \operatorname{GL}_{d+1}(\mathbb{Q}^{\operatorname{ab}})$ such that $\bar{\xi}(\sigma) = \sigma(A)A^{-1}$ for every $\sigma \in \operatorname{Gal}(\mathbb{Q}^{\operatorname{ab}}/\QQ)$.

\indent Let $f \in M_{k,G}$ be a modular form. Let $v$ be its corresponding coordinate vector with respect to the basis $\mathcal{B}$. Consider the modular form $g\in M_{k,G}\otimes_\QQ \QQ^{\operatorname{ab}}$ with the coordinate vector $ v \cdot A^{-1} \in (\mathbb{Q}^{\operatorname{ab}})^{d+1}$. The $\bullet$ action on $M_{k,G}\otimes_\QQ \QQ^{\operatorname{ab}}$ gives an action on the vector space $(\mathbb{Q}^{\operatorname{ab}})^{d+1}$ with respect to the basis $\mathcal{B}$. For all $\sigma \in \operatorname{Gal}(\mathbb{Q}^{\operatorname{ab}}/\QQ)$, we have that
\[
\sigma \bullet (v\cdot A^{-1})=\sigma(v\cdot A^{-1})\cdot\bar{\xi}_\sigma=v\cdot\sigma(A^{-1})\cdot\bar{\xi}_\sigma=v\cdot\sigma(A)^{-1}\sigma(A)A^{-1}=v\cdot A^{-1}.
\]

Hence, the modular form $f_\gamma\colonequals g$ with the coordinate vector $v \cdot A^{-1}$ is a modular form in $M_{k,G_\gamma}$.
Applying the matrix $A^{-1}$ to the basis column $\mathbf f$, we obtain $\mathbf f'=A^{-1}\mathbf f$, whose entries form a basis $\mathcal{B}'$ of $M_{k,G_\gamma}$.

\subsection{Twisting the models of modular curves} 
Assume that we have an explicit smooth projective model $C\subseteq \mathbb{P}_{\QQ}^r$ for $X_G$, obtained from a very ample sheaf as in \S \ref{S:Models},  where $G\in \scrF(\calG,B)$. The model $C$ arises from linearly independent modular forms $f_0,\ldots,f_r\in V\subseteq M_{k,G}$ as explained in \S \ref{S:Models}. The $\QQ$-vector space $V$ is chosen so that there is an action of $\calG/G$ on $V$. In particular $C$ is defined by $F_1,\ldots, F_s \in \QQ[x_0,\ldots , x_r]$ where $F_i(f_0,\ldots , f_r)=0$. Let $\gamma : \Zhat^{\times} \rightarrow \calG/G$ be a continuous homomorphism and let $\xi:=\gamma \circ \chi_{\operatorname{cyc}}$. Let $(F_i)_\xi:=F_i((x_0,\ldots,x_r)\cdot A^T)$. 

\begin{thm}\label{T:twist}
Let $I$ be the ideal $I\subseteq \QQ^{\operatorname{ab}}[x_0,\ldots , x_r]$ generated by $(F_1)_\xi,\ldots , (F_s)_\xi$. Let $C'$ be the curve given by the zero set of $I$. The ideal $I$ is defined over $\QQ$ and the curve $C'$ is isomorphic to the twist of $C$ by $\xi$. In particular, $C'$ is a projective model of $X_{G_\gamma}$.
\end{thm}
\begin{proof}
First, observe that $\xi_\sigma$ is an automorphism of $X_G$ defined over $\QQ$ and so each $\overline{\xi}(\sigma)$ induces an automorphism of the model $C$ of the modular curve $X_G$. The map $\overline{\xi}$ is a group homomorphism and substitution by $\overline{\xi}(\sigma)$ preserves the homogeneous ideal $I(C)$.

Writing $\mathbf x=(x_0,\ldots,x_r)^{\mathsf T}$, we have
\[
\sigma(F_i(A\mathbf x))=F_i(\bar\xi(\sigma)A\mathbf x).
\]
Thus $I$ is Galois-stable and descends to $\QQ$: if $I_0:=I\cap\QQ[x_0,\ldots,x_r]$, then $I=I_0\QQ^{\operatorname{ab}}[x_0,\ldots,x_r]$. The individual generators $(F_i)_\xi$ need not have rational coefficients.

Let $\mathcal{B}:=\{f_0,\ldots,f_r\}$ be the basis of $V$ as above. Applying $A^{-1}$ to the basis column as in \ref{SS:getting a basis}, we get a basis $\mathcal{B}'$ for a vector space $V_\gamma$ which is acted on by $\calG$. The space $V_\gamma$ is associated to a sheaf $\calF_\gamma$ on $X_{G_\gamma}$ which is very ample. The polynomials $(F_i)_{\xi}$ vanish on the basis $\mathcal{B}'$, so the descended curve $C'$ is isomorphic to $X_{G_\gamma}$ over $\QQ$.
Over $\QQ^{\operatorname{ab}}$, the isomorphism $\varphi:C\rightarrow C'$, $x\mapsto A^{-1}x$, satisfies
\[
\varphi^{-1}\circ\sigma(\varphi):x\longmapsto A\sigma(A)^{-1}x=\bar\xi(\sigma)^{-1}x,
\]
which is the action of $\xi_\sigma$ on points under our identification with curve automorphisms.
\end{proof}
\subsubsection{Computing the matrix $A$}\label{H90matrix} In practice, we work over the field $\QQ(\zeta_N)$ to obtain a basis for $M_{k,G_\gamma}$ where $N$ is the least common multiple of the level of $G$ and the conductor of $\gamma$ (the least positive integer $c$ such that $\gamma$ factors through $(\ZZ/c\ZZ)^{\times}$). The Hilbert 90 Theorem states that $H^1(\Gal(\QQ(\zeta_N)/\QQ),\operatorname{GL}_n(\QQ(\zeta_N)))$ is trivial so given the cocycle $\eta(\sigma):=\bar\xi(\sigma)^{\mathsf T}$, there exists a matrix $B\in \operatorname{GL}_n(\QQ(\zeta_N))$ such that $\eta(\sigma)=B^{-1}\sigma(B)$. We set $A=B^{\mathsf T}$, so that $\bar\xi(\sigma)=\sigma(A)A^{-1}$. \\
\indent The matrix $B$ can be explicitly computed. In practice, we are using the algorithm and implementation given in \cite{rakvi}*{\S 5.3}.

\section{The Algorithm} \label{S:algorithm}
\indent In this section, we put together the work done in previous sections and describe an algorithm to compute a projective model of a modular curve $X_G$ where $G\subseteq \Ggl{\Zhat}$ is an open subgroup with genus at most $g$ for a fixed natural number $g$. \\
\indent Our algorithm has two parts. The families mentioned in Theorem \ref{T:main familiesend} must be computed along with a chosen representative for each family. Using these precomputed data, we then provide an algorithm that, given a modular curve $X_G$, finds the family of groups it lies in, computes the cocycle with respect to the representative, and twists the representative curve to get a projective model of $X_G$. 

\begin{algorithm}[\textbf{Precomputation}]\label{alg1}
\begin{romanenum}
\item Compute $\scrA_g$, a set of representatives up to conjugacy in $\GL_2(\Zhat)$ of the agreeable subgroups of $\Ggl{\Zhat}$ whose genus is at most $g$. In subsection \ref{SS:computeagg}, we have described an algorithm for computing them.
\item  For each $\calG \in \mathscr{A}_g$, compute the subgroups $B$ such that $[\calG,\calG]\subseteq B \subseteq \calG \cap \Ssl{\Zhat}$. Note that $[\calG,\calG]$ and $\calG \cap \Ssl{\Zhat}$ are open subgroups of $\operatorname{SL}_2(\Zhat)$, so for each agreeable subgroup there are only finitely many such subgroups $B$. Retain only those $B$ with $-I\in B$ and $B\cap\SL_2(\ZZ)$ of genus at most $g$.
\item Form all the possible families (up to conjugacy) $\scrF(\calG,B)$ from the retained pairs.
\item For each such family, determine if it is empty. If it is not empty find a representative $W\in \scrF(\calG,B)$. Discard the empty families; the remaining collection is $\scrF_g$.
\item For each family $\scrF(\calG,B)\in \scrF_g$ with a representative $W$, compute a model $C\subseteq \PP_\QQ^d$ of $X_W$ and the associated relative $j$-map $\pi_{\calG/W}:X_W \rightarrow X_\calG$ via the methods described in \S \ref{S:Background} and the algorithm given by Zywina in \cite{explicitopenimage}. In particular, we compute a set of linearly independent modular forms $f_0,...,f_d \in M_{k,W}$ and $C$ is defined by homogeneous polynomials $F_1,\ldots , F_s \in \QQ[x_0, \ldots , x_d]$ such that $F_i(f_0,\ldots , f_d)=0$ for $i=1,\ldots,s$. 

\end{romanenum}
\end{algorithm}
\begin{remark} 
\leavevmode
\begin{itemize}
\item If two open subgroups $G,G'\subseteq \GL_2(\Zhat)$ are conjugate to each other in $GL_2(\Zhat)$, then $X_G$ and $X_{G'}$ are isomorphic over $\QQ$. Therefore, it is enough to consider the agreeable subgroups and families of modular curves up to conjugacy in $\GL_2(\Zhat)$.
\item For the first step, we start from the data of congruence subgroups of $\operatorname{SL}_2(\ZZ)$ which is given in \cite{MR2016709} for $g\le 24$. Hence we obtain the set $\scrA_g$ for $g=24$. 
\item Theorem \ref{T:representativeinfam} gives a criterion for checking whether a family is empty or not. Based on this theorem, an implementation is given in \cite{openimagenumberfields} to find a representative in $\scrF(\calG,B)$.
\item Algorithm \ref{alg1}(v) is, computationally, the most expensive part of the precomputation as it includes the computation of Eisenstein series (and their $q$-expansions for possibly high precision) that span $M_{k,W}$ for a certain $k\in \NN$.
\end{itemize}
\end{remark}

\noindent \textbf{Computing the models:} After the precomputation, one can use the following algorithm to compute a projective model of $X_G$. It takes an open subgroup $G\subseteq \operatorname{GL}_2(\Zhat)$ with full determinant and $-I\in G$ as input.

\begin{algorithm}\label{mainalgo}
Let $G\subseteq \Ggl{\Zhat}$ be an open subgroup where $X_G$ has genus at most $g$. This algorithm computes a model $C\subseteq \PP_{\QQ}^d$ of the modular curve $X_G$.
\begin{romanenum}
\item Compute $B:=G \cap \Ssl{\Zhat}$ and the agreeable closure of $G$, which we call $\calG$.
\item By Theorem \ref{T:main familiesgroups}, $\calG$ is in $\scrA_g$, possibly up to conjugacy in $\GL_2(\Zhat)$. If necessary, conjugate $G$, $\calG$ and $B$, and replace them with their suitable conjugates so $\calG\in \scrA_g$. Find the family $\scrF(\calG,B)\in \scrF_g$.	
\item Since $\scrF(\calG,B)$ is not empty, we have precomputed a representative $W\in \scrF(\calG,B)$. Compute the homomorphism $\gamma: \Zhat^{\times} \cong G/B \rightarrow \calG/W$. Then $G$ is equal to $W_\gamma$ by Theorem \ref{familygroups}.
\item Compute the associated cocycle $\xi$ by precomposing with the cyclotomic character $\chi_{\operatorname{cyc}}$ and the related cocycle $\bar{\xi}$ as in \S \ref{SS:getting a basis}. 
\item Compute the Hilbert 90 matrix $A$, as described in \S\ref{H90matrix}, over a finite cyclotomic field $K$.
\item Apply the matrix $A^T$ to the polynomials defining $X_W$ to get polynomials $(F_1)_{\xi},\cdots , (F_s)_{\xi}$. These have coefficients in $K$, and their generated ideal descends to $\QQ$.
\item By Theorem \ref{T:twist}, the descended curve $C_\xi\subseteq\PP^d_\QQ$ is a projective model of $X_G$. By \cite{explicitopenimage}*{\S\S 5.3.1--5.3.2}, $I(C_\xi)$ is generated by $I(C_\xi)_2 \cup I(C_\xi)_3 \cup I(C_\xi)_4$. For each $m\in\{2,3,4\}$, find a $K$-basis $g_1,\ldots,g_t$ of the full degree-$m$ part of the transformed ideal, including suitable monomial multiples of lower-degree generators. For a $\QQ$-basis $\theta_1,\ldots,\theta_e$ of $K$, the coefficientwise traces $\operatorname{Tr}_{K/\QQ}(\theta_a g_j)$ span $I(C_\xi)_m$. Select a $\QQ$-basis from these traces by linear algebra.
\end{romanenum}
\end{algorithm}

\begin{remark}
 Let $G$ be the input of our algorithm, and let $N$ be its level. Let $N_1,N_2$ be the levels of $\calG$ and $W$, respectively. The level of $G$ is not bounded in terms of $N_1$ and $N_2$, it can be arbitrarily big. We do most of our group computations modulo $\operatorname{lcm}(N_1,N_2)$. To compute $\gamma$, we can use the common modulus $L=\operatorname{lcm}(N,N_2)$, so that $\gamma:(\ZZ/L\ZZ)^{\times}\rightarrow\calG/W$. The splitting field of the cocycle, and hence the Hilbert 90 and trace computations, can also depend on the input level $N$.
\end{remark}
This algorithm has been implemented for $g\leq12$ in \cite{repo}. Combining our algorithm with previous results, we restate Theorem \ref{T:main families1}.

\begin{thm}\label{T:main familiesend}
Fix a non-negative integer $g$. Let $G$ be an open subgroup of $\GL_2(\Zhat)$ of genus $g$ with full determinant and $-I\in G$.
\begin{enumerate}
    \item There are only finitely many families of modular curves of genus $g$. These families are effectively computable.
 \item Let $G$ be an open subgroup of $\GL_2(\Zhat)$ of genus $g$ with full determinant and $-I\in G$. There is an effective algorithm utilizing families of modular curves, which takes as input the group $G$ and outputs a projective curve $C_{G}\subseteq \PP_\QQ^r$ for some $r> 0$ such that $C_{G}$ is isomorphic to $X_{G}$. 
\end{enumerate}

\end{thm}
\begin{proof}
 By Corollary \ref{C:twist}, the family of groups mentioned in Theorem \ref{T:main familiesgroups} corresponds to a family of twists of modular curves which proves (1). The algorithm is described in Algorithm \ref{mainalgo}.

\end{proof}

\appendix

\section{Computing Models of Universal Elliptic Curves (with Yongyuan Huang and Rakvi)}

From the moduli-theoretic viewpoint, $X_G$ is in general only a coarse moduli space. In particular, if $-I\in G$, the stabilizer of every noncuspidal geometric point contains $\ZZ/2\ZZ$. However, when $-I\notin G$ and the level of $G$ is greater than or equal to $3$, away from the cusps and elliptic points the moduli problem is fine, and hence there exists a \textit{universal elliptic curve} over $U=X_G-\pi_{G}^{-1}(\{0,1728,\infty\})\subseteq X_G$, where $\pi_G:X_G\to\PP^1$ is the $j$-map, in the following sense.

\begin{defn} Following \cite[\S13]{explicitopenimage}, we say that an elliptic scheme $\mathcal{E}_G=\mathcal{E}\to U$ is a \textit{universal elliptic curve} over $U$ if the following hold for any number field $K$:
\begin{enumerate}
    \item for any point $t\in U(K)$, $E_t$, the fiber of $\mathcal{E}\to U$ over $t$, has $j$-invariant $\pi_G(t)$,
    \item for all $t\in U(K)$, $\rho_{E_t}^{\ast}(\Gal_K)$ is conjugate in $\GL_2(\Zhat)$ to a subgroup of $G$.
\end{enumerate}
\end{defn}

To further elucidate the distinction, suppose $E/K$ is a non-CM elliptic curve such that $j(E)\in\pi_G(X_G(K))$. When $-I\in G$, the property $\rho_E^{\ast}(\Gal_K)\subseteq G$ is preserved by quadratic twisting of $E$. When $-I\notin G$, fix compatible torsion bases for which $\rho_E^{\ast}(\Gal_K)\subseteq\langle G,-I\rangle$. There is then a unique quadratic twist $E_d$ of $E$ such that $\rho_{E_d}^{\ast}(\Gal_K)\subseteq G$ with respect to the transported bases; cf. \cite[Lemma 5.1 and Remark 5.2]{rousezureick}.

From now on, we suppose that $G\subseteq\GL_2(\Zhat)$ is open, $\det(G)=\Zhat^{\times}$, and $-I\notin G$; its level is then at least $3$. Let $f_0\in M_{3,G}$ be a nonzero weight $3$ modular form, which exists by \cite[Lemma 4.6]{explicitopenimage}. Let $E_6$ be the normalized Eisenstein series of weight $6$, and let $j$ be the $j$-map as usual. Let $0\neq r\coloneq f_0^2/E_6\in\QQ(X_G)$. The equation
\begin{equation}\label{unimodel}
    y^2=x^3-27\cdot r^2\left(1-\frac{1728}{j}\right)x+54\cdot r^3\left(1-\frac{1728}{j}\right)^2
\end{equation}
gives a model for the generic fiber of a universal elliptic curve over $U$; after suitable local changes of coordinates, it extends to an elliptic scheme over $U$ \cite[Prop. 13.2]{explicitopenimage}. Thus in order to compute a model for $\mathcal{E}_G$, it suffices to explicitly compute $r$ and $j$ as elements of $\QQ(X_G)$. In the next section, we explain how to compute a presentation of an arbitrary quotient of modular forms of the same weight.

\subsection{Computing rational functions on modular curves}

Let $G\subseteq \GL_2(\Zhat)$ be an open subgroup satisfying $\det(G)=\Zhat^{\times}$ and $-I \notin G$. Let $G^+\colonequals \langle G,\pm I\rangle$. Let $M_{k,G}$ be the space of modular forms of weight $k$ defined in \S\ref{SS:MkG}. Let $u,v \in M_{k,G}$ be nonzero, and let $\eta\colonequals \frac{u}{v}$ denote their quotient. The case $u=0$, $v\ne0$, is represented by $0/1$.

We assume that a smooth projective model has been computed for $X=X_G=X_{G^+}$, and let $f_1,\ldots,f_n$ be a $\QQ$-basis of $H^0(X,\calF)$, where $\calF$ is the line bundle used to obtain the model, cf. \S\ref{S:Models}. We call them the \defi{model sections} and write $f=(f_1,\ldots,f_n)$. The space $H^0(X,\calF)$ is identified with a subspace of $M_{\kappa,G^+}$, where $\kappa$ is the common even weight of the modular forms $f_i$. We assume that the evaluation map in \eqref{E:normalgeneration} is an isomorphism, so that every section of $\calF^{\otimes d}$ is represented by a homogeneous degree $d$ polynomial in the model sections. This holds for both the nonhyperelliptic canonical models and the models with $\deg\calF\geq2g+1$ described in \S\ref{L:computing a model}.

Since $u$ and $v$ have the same weight and both belong to $M_{k,G}$, their quotient $\eta$ is a rational function on $X_G$. It follows that there exist an integer $d\geq0$ and homogeneous polynomials $P,Q\in\QQ[x_1,\ldots,x_n]_d$, with $Q(f)\ne0$, such that
\begin{equation}\label{Eq:quotient}
    \ddfrac{P(f_1,\ldots,f_n)}{Q(f_1,\ldots,f_n)}=\eta.
\end{equation}

The integer $d$ in \eqref{Eq:quotient} is called a \defi{presentation degree}. The goal of this section is to give a computable degree at which such a presentation of $\eta$ is guaranteed to exist. In Section \ref{Appe:algo}, we will describe an algorithm for computing a presentation of degree $d$. The algorithm requires explicit computations with the Fourier expansions of $u,v$ and the model sections. One needs to find an integer $p$, called a \defi{precision bound}, for which the truncated expansions are sufficient to establish exact equalities. We tackle this problem in Section \ref{Appe:precision}.

\begin{remark}
    More generally, if $\eta\in\QQ(X_G)$ is a rational function on $X_G$, then it can be written as a quotient of the form \eqref{Eq:quotient}. In particular, we can describe the $j$-map $\pi_G:X_G\rightarrow\PP^1$ as a quotient of two polynomials of the same degree evaluated at the model sections. The algorithm described in the next section can be used to compute the $j$-map by taking $u=E_4^3$ and $v=\Delta$, where $\Delta=(E_4^3-E_6^2)/1728$.
\end{remark}
\subsubsection{Notation:} Let $\mu=[\SL_2(\ZZ): \pm \Gamma_G]$, and let $m$ be the level of $\Gamma_G$. Let $\mathcal{L}$ be a line bundle on $X$, where $X$ is a smooth, projective, geometrically integral curve. In our work, $X=X_G$ and we write $g$ for its genus. Let $e$ be a nowhere zero local regular section of $\calL$, which is called a \defi{frame}. Every local regular section on the same open set can be written uniquely in the form $ae$, where $a$ is a regular function on that open set.

For a nonzero global regular section $s$ of $\calL$, which has no poles, denote by $Z_s$ its zero divisor. We have $\mathcal{O}_X(Z_s)\simeq \mathcal{L}$. The degree of $\calL$ is equal to the degree of $Z_s$ as a divisor. For a divisor $A$ on the curve, we define $\mathcal{L}(A)\colonequals \mathcal{L} \otimes \mathcal{O}_X(A)$. Global regular sections of $\calL(A)$ are meromorphic sections $s$ of $\mathcal{L}$ with $\operatorname{div}(s)+A \geq 0$. If $A$ is effective, then $\calL (-A)$ simply imposes vanishing along $A$. The model degree $D>0$ of the curve $X_{G}$ is defined as the degree of the line bundle $\calF$.

For an even integer $a\geq0$, define $\scrL_a=\Omega_X^{\otimes a/2}(\calD_a)$, using the algebraic line-bundle notation of \S\ref{S:Models}. Here $\calD_a$ is the divisor in \eqref{E:divisor}, with coefficient $a/2$ at each cusp and $\lfloor (a/2)(1-1/e)\rfloor$ at each elliptic point of order $e$. Define
\begin{equation}\label{E:linebundledegree}
    B_a\colonequals \deg \scrL_a = \frac{a}{2}(2g-2+c)+v_2\lfloor a/4\rfloor +v_3 \lfloor a/3\rfloor,
\end{equation}
where $c$ is the number of geometric cusps, and $v_2$ and $v_3$ are the numbers of elliptic points of orders $2$ and $3$, respectively. The cusp count denoted by $r$ in \S\ref{S:Models} is denoted by $c$ here, since $r=f_0^2/E_6$ is already in use.

\subsubsection{Finding the degree}
We return to the discussion of rational functions. Let $H=\operatorname{div}_{\infty}(\eta)$ be the pole divisor of $\eta$ and let $p_r\colonequals\deg H$. For even $k$, view $u,v$ as global sections of $\scrL_k$, let $Z_u,Z_v$ be their zero divisors, and define $C_{u,v}=\min(Z_u,Z_v)$ coefficientwise. One has
\begin{equation}\label{Eq:degreebound}
    H=Z_v - C_{u,v}, \quad \quad p_r =\deg \scrL_k - \deg C_{u,v} \le \deg \scrL_k.
\end{equation}

Define
\begin{equation}\label{eq:Sk}
 S_k=
 \begin{cases}
 B_k,&k\ \text{even},\\
 B_{2k}/2,&k\ \text{odd}.
 \end{cases}
\end{equation}

We can bound the degree of the pole divisor of $\eta$ in terms of $S_k$.
\begin{prop}
    For nonzero $u,v \in M_{k,G}$,
    \[
    p_r\leq B\colonequals \lfloor S_k \rfloor \leq \lfloor \ddfrac{k\mu}{12}\rfloor.
    \]
\end{prop}

\begin{proof}
    If $k$ is even, equation \eqref{Eq:degreebound} implies the first inequality. For odd $k$, one has $u^2,v^2\in M_{2k,G^+}$, hence
    \[
    2p_r=\deg (\operatorname{div}_{\infty}(\eta^2))\leq B_{2k}.
    \]
    For the second inequality, observe that the genus formula for $X_{G}$ gives 
    \[
    \ddfrac{a\mu}{12}-B_a=v_2 \left\{ a/4\right\} +v_3 \left\{ a/3\right\} \geq 0
    \]
    where $a$ is an even integer and the curly brackets denote the fractional part. Applying this with $a=k$ for even $k$ and $a=2k$ for odd $k$ proves the assertion.
\end{proof}

The next proposition provides an upper bound on the presentation degree.

\begin{prop}\label{prop:degree}
Let $d\geq0$. A degree $d$ presentation of $\eta$ exists if and only if $H^0(X,\calF^{\otimes d}(-H))\ne0$. A presentation exists by degree
\[
d=\max \{1,\lceil\ddfrac{B+g}{D}\rceil\}.
\]
\end{prop}
\begin{proof}
    Write $P(f)/Q(f)=\eta$, and let $q=Q(f)\ne0$. The section $q$ vanishes at every pole of $\eta$ to at least that pole order. Hence $q\in H^0(X,\calF^{\otimes d}(-H))$. Conversely, if there exists such a nonzero section $q$, then $q\eta$ is a regular section in $H^0(X,\calF^{\otimes d})$. The isomorphism \eqref{E:normalgeneration} provides the polynomials $P$ and $Q$.

    The degree of $\calF^{\otimes d}(-H)$ is $dD-p_r$. It is necessary that $dD\geq p_r$ for a presentation of degree $d$ to exist. By Riemann--Roch, $h^0(X,\calF^{\otimes d}(-H))\geq dD-p_r+1-g$, so it suffices to have $dD-p_r\geq g$. For the last statement, use $p_r\leq B$.
\end{proof}

\begin{remark}\label{rem:checking}
    Assume that we have a candidate rational function $a=P(f)/Q(f)$, with $P,Q$ homogeneous of degree $d$ and $Q(f)\ne0$, and want to check if $a=\eta$. If $a\ne\eta$, we have
\[
\operatorname{div}_\infty(a-\eta)\leq\operatorname{div}_\infty(a)+\operatorname{div}_\infty(\eta),
\]
which implies $\deg\operatorname{div}_\infty(a-\eta)\leq dD+B$. A nonzero rational function on a smooth projective curve has zero and pole divisors of equal degree. Thus, if $a-\eta$ has zeros of total degree greater than $dD+B$, counted with multiplicity, then $a=\eta$.
\end{remark}

\subsection{Algorithm for computing rational functions}\label{Appe:algo}
The previous section proved that a presentation of $\eta$ exists by a computable degree $d$. In this section, we describe how to find such a presentation.

Fix the curve model and its sections $f_1,\ldots,f_n$. Fix $u,v$ and $\eta=u/v\ne0$ as above, and a positive integer $d$. Write $f=(f_1,\ldots,f_n)$. Let $I=I(X_G)$ be the full homogeneous ideal defining the chosen model of $X_G$. Let $R\colonequals\QQ[x_1,\ldots,x_n]/I$ and let $R_d$ denote its degree $d$ component. The polynomials $P$ and $Q$ defining a degree $d$ presentation determine classes in $R_d$.

Consider the linear map $\calR:R_d\times R_d\rightarrow M_{\ell,G}$, where $\ell=k+d\kappa$, defined by
\[
\mathcal{R}(P,Q)=vP(f)-uQ(f).
\]
Its kernel is
\[
\mathcal{K}_d\colonequals\{(P,Q)\in R_d\times R_d\colon vP(f)-uQ(f)=0\}.
\]
The proof of Proposition \ref{prop:degree} shows that, under the evaluation isomorphism, the assignment $q\mapsto(\eta q,q)$ identifies $H^0(X_G,\calF^{\otimes d}(-H))$ with $\mathcal{K}_d$.

Since $u,v\ne0$, every nonzero pair $(P,Q)\in\mathcal{K}_d$ has $Q(f)\ne0$ and gives a presentation $P(f)/Q(f)=\eta$. Our goal is therefore to explicitly find a nonzero element of $\mathcal{K}_d$.

Let $h_d=\dim_{\QQ}R_d$. Choose degree $d$ monomials $b_1,\ldots,b_{h_d}$ whose classes form a basis of $R_d$. Note that evaluating $b_\alpha$ at the model sections gives a modular form $b_\alpha(f)=b_\alpha(f_1,\ldots,f_n)$ of weight $d\kappa$.

Using these monomials as representatives, write
\[
P=\sum_{i=1}^{h_d}a_ib_i\quad \text{and} \quad Q=\sum_{i=1}^{h_d}c_ib_i,
\]
where $(a_i)_{i=1,\ldots,h_d}$ and $(c_i)_{i=1,\ldots,h_d}$ are $2h_d$ unknown rational coefficients. The equation $\mathcal{R}(P,Q)=0$ is linear in these coefficients. We determine the resulting linear equations using the Fourier expansions at the cusps.

Fix a cusp $P_i$, a representative $A_i\in\SL_2(\ZZ)$ with $A_i\infty=P_i$, and the parameter $s_i=q_{w_i}$, where $w_i$ is the cusp width for $\Gamma_G$. Write
\[
u_i=u\vert_k A_i, \quad v_i=v\vert_k A_i, \quad \beta_{\alpha,i}=b_{\alpha}(f) \vert_{d\kappa} A_i.
\]
These are power series in the same parameter $s_i$. We have
\[
\mathcal{R}(P,Q)\vert_{\ell}A_i=\sum_{\alpha=1}^{h_d}a_\alpha v_i \beta_{\alpha,i}-\sum_{\alpha=1}^{h_d}c_\alpha u_i \beta_{\alpha,i}.
\]
In order to obtain a linear system of equations, we use the modular forms
\begin{equation}\label{eq:modformcolumns}
    v_i\beta_{1,i},\ldots,v_i\beta_{h_d,i}, -u_i\beta_{1,i},\ldots,-u_i\beta_{h_d,i}
\end{equation}
as columns to expand. Choose an integer cutoff $N_i$ and assume that all these power series are known modulo $s_i^{N_i}$, that is, at exponents less than $N_i$. For every $0\leq n<N_i$, we impose
\[
\sum_\alpha a_\alpha\coef_n(v_i\beta_{\alpha,i}) -\sum_\alpha c_\alpha\coef_n(u_i\beta_{\alpha,i})=0.
\]
Here $\coef_n$ extracts the coefficient of $s_i^n$. The coefficients initially lie in $\QQ(\zeta_N)$, where $N$ is the level of $G$. Expand every coefficient in a fixed $\QQ$-basis $\xi_1,\ldots,\xi_{\varphi(N)}$ of this field and equate each rational coordinate to zero. Combining these equations at each cusp gives a rational matrix $A_d((N_i)_i)$. Its $2h_d$ columns correspond to the products in \eqref{eq:modformcolumns}.

The rows of $A_d((N_i)_i)$ are indexed by three parameters: the cusp index $i$, the exponent $n$, and the basis index $j\in\{1,\ldots,\varphi(N)\}$. For
\[
x=(a_1,\ldots,a_{h_d},c_1,\ldots,c_{h_d})^{\intercal},
\]
the equation $A_d((N_i)_i)x=0$ means exactly that $\mathcal{R}(P,Q)$ vanishes to the prescribed precision. In general, we only have an inclusion $\mathcal{K}_d\subseteq\ker A_d((N_i)_i)$, but this can be made to be an equality by choosing the precision carefully. A nonzero vector $x$ satisfying $A_d((N_i)_i)x=0$ and $Q(f)\ne0$ gives a candidate for $\eta$ in the form $P(f)/Q(f)$.

Hence, we choose such a vector and form the candidate function for $u/v$. One can then check whether it equals $\eta$ by Remark \ref{rem:checking}. 

\begin{remark}\label{rem:minmaxdegree}
    For $f_0^2/E_6$, the implementation starts testing for a presentation at $d=\max\{1,\lceil B_6/D\rceil\}$ and increases the degree and precision as needed. This starting degree need not be minimal: the necessary condition in Proposition \ref{prop:degree} is $dD\geq p_r$, whereas $B_6$ is an upper bound on $p_r$.
\end{remark}

\subsection{Precision}\label{Appe:precision}

We express modular forms using their power series expansions at the cusps. If $s$ is the parameter of an expansion, then the expression
\[
f=\sum_{n=0}^{L-1}a_ns^n+O(s^L)
\]
specifies the coefficients with exponents less than $L$. The goal of this section is to find a positive integer $p(d)$ for which the truncated expansions determine the exact relation space $\mathcal{K}_d$, and hence guarantee the correctness of the output of the preceding algorithm.

Let $m$ be the level of $\Gamma_G$, and let $w_i$ and $h_i$ be the widths of the cusp $P_i$ for $\Gamma_G$ and $\pm\Gamma_G$, respectively. Put $\epsilon_i=w_i/h_i\in\{1,2\}$. A local parameter on the coarse curve is $t_i=s_i^{\epsilon_i}$. For an integer $p\geq1$, set $N_i(p)=\lceil w_ip/m\rceil$, determining the cutoff at each cusp from a single parameter. For a nonzero modular form $\Phi\in M_{\ell,G}$, define its normalized cusp order by 
\[
\nu_i(\Phi)=\ddfrac{\ord_{s_i}(\Phi\vert_{\ell} A_i)}{\epsilon_i}.
\]
The sum of these orders is bounded by
\[
\sum_{i=1}^c \nu_i(\Phi)\leq S_{\ell}.
\]
For even $\ell$, this follows from the degree of $\scrL_{\ell}$. For odd $\ell$, apply the even-weight bound to $\Phi^2$ and divide by $2$. This bound can be used to obtain the required precision.

\begin{thm}
    Fix a positive integer $d$ and a basis of $R_d$, and put $\ell=k+d\kappa$. Suppose that at each cusp $P_i$, the expansions of $u,v$ and every model section $f_j$, expressed in the parameter $s_i=q^{1/w_i}$, are known modulo $s_i^{N_i(p)}$. If 
    \[
    \sum_i \ddfrac{N_i(p)}{\epsilon_i}> S_{k+d\kappa}
    \]
    then the coefficient matrix $A_d((N_i(p))_i)$ has kernel equal to $\mathcal{K}_d$. In particular,
    \[
    p(d)=\left\lfloor\frac{mS_{k+d\kappa}}{\mu}\right\rfloor+1
    \]
    is sufficient to recover the exact relation space $\mathcal{K}_d$ from these expansions.
\end{thm}
\begin{proof}
    Since the input forms are holomorphic, their products are known modulo $s_i^{N_i(p)}$. Let $x$ be a vector in the kernel of $A_d((N_i(p))_i)$, and let $P,Q$ be the corresponding polynomials. Then $\nu_i(\mathcal{R}(P,Q))\geq N_i(p)/\epsilon_i$ at each cusp, with the convention that the zero form has infinite order. 

    Assume $\mathcal{R}(P,Q)$ is nonzero. Summing over all cusps gives
    \[
        S_{\ell}<\sum_i\frac{N_i(p)}{\epsilon_i}\leq\sum_{i=1}^c\nu_i(\mathcal{R}(P,Q))\leq S_{\ell},
    \]
    a contradiction. Conversely, every exact relation has all these coefficients equal to zero, so the two kernels agree.

    To obtain the sufficient precision, note that $\sum_iw_i/\epsilon_i=\mu$ and
    \[
    \sum_i \ddfrac{N_i(p)}{\epsilon_i}\geq \ddfrac{p}{m}\sum_i w_i/\epsilon_i=\ddfrac{p\mu}{m}.
    \] The stated value of $p(d)$ makes the last quantity strictly greater than $S_{k+d\kappa}$.
\end{proof}

We can also give a sufficient precision using only the level and weights.
\begin{cor}

    Let $\kappa$ be the weight of the modular forms defining the model. Let $d$ be the degree at which we are searching for relations, and let $k$ be the weight of the modular forms whose ratio we are computing. Under the hypotheses of the preceding theorem, it suffices to take
    \[
    \lfloor \ddfrac{m(k+d\kappa)}{12}\rfloor+1
    \]
    as the precision integer, where $m$ is the level of $\Gamma_G$.
\end{cor}
\begin{proof} This follows from $S_{k+d\kappa}\leq(k+d\kappa)\mu/12$. \end{proof}

\subsection{The algorithm for computing the universal elliptic curve}

Here we explain our algorithm for computing a model over $\QQ(X_G)$ for the generic fiber of the universal elliptic curve.

\begin{enumerate}
    \item As input, we take an open subgroup $G$ of $\GL_2(\Zhat)$ of full determinant that does not contain $-I$, a smooth projective model for the corresponding coarse curve, its $j$-map, and $r=f_0^2/E_6$. The model and both functions are given with respect to the same identification with $X_G$.

    \item We then construct the initial model $y^2=x^3+Ax+B$, where \[A=-27r^2\big(1-\frac{1728}{j}\big),\qquad B=54r^3\big(1-\frac{1728}{j}\big)^2.\] In \texttt{FirstReduction} we implement an algorithm to find a suitable $u\in\QQ(X_G)^{\times}$ that simplifies the coefficients $A/u^4$ and $B/u^6$. It does this by constructing the divisor \[\sum_{P}\min\bigg\{\bigg\lfloor\frac{\nu_P(A)}{4}\bigg\rfloor,\bigg\lfloor\frac{\nu_P(B)}{6}\bigg\rfloor\bigg\}P,\] where $P$ runs over the closed points of $X_G$ and $\nu_P$ is the order at $P$. It then applies divisor reduction relative to a base divisor of positive degree obtained from \texttt{LowDegreeDivisor}, using \texttt{Magma}'s \texttt{Reduction(D,E)} routine \cite{magma-divisor-reduction}, with $D$ the displayed divisor and $E$ the base divisor. The routine returns $\widetilde D,n,E,a$ satisfying $D=\widetilde D+nE-\operatorname{div}(a)$. Taking $u=a^{-1}$ gives the corresponding change of Weierstrass equation.

    \item The function \texttt{torsionorder} determines the largest order $d$ of a rational torsion point on the generic fiber, returning $1$ if the rational torsion is trivial. The function \texttt{FinalReduction} chooses the model depending on $d$. If $d=1$, the final Weierstrass equation has the form $y^2=x^3+Ax+B$. If $d=2$, it has the form $y^2=x^3+a_2x^2+a_4x$. Otherwise, it has the form $y^2+a_1xy+a_3y=x^3+a_2x^2$.

    \item Finally, the algorithm checks that the resulting model is nonsingular, that it has the correct $j$-invariant, that its $c_4$ and $c_6$ invariants differ from those of the initial model by the appropriate fourth and sixth powers of a common nonzero scaling factor, and that any distinguished point has the asserted exact torsion order.
    
\end{enumerate}

\subsection{Example} We will conclude with an explicit example. Consider the group $G$ whose associated modular curve has LMFDB label \texttt{30.48.1-30.i.2.4}. This is a genus $1$ curve with infinitely many rational points. A model for the generic fiber of its universal elliptic curve is \[Y^2+a_1(s,t)XY+a_3(s,t)Y=X^3+a_2(s,t)X^2,\] where

\begin{align*}
    a_1(s,t) &=\frac{120s^2}{t^2} + \frac{435s}{t} + \frac{300t^2-64}{t^2}, \\
    a_2(s,t) &= \frac{-1800s^2}{t^2} - \frac{6525s}{t} + \frac{-4725t^2 + 960}{t^2}~\text{and}\\
    a_3(s,t) &= \frac{27000s^2}{t^2} + \frac{97875s}{t} + \frac{70875t^2 - 14400}{t^2}.
\end{align*}

Further, a coarse model in $\PP^2$ is given by \[C\colon-15a^3-60a^2b-60ab^2+8ac^2-15b^3+3bc^2=0,\] where $s=a/c$ and $t=b/c$. The point $[3/10:-7/10:1]$ is a rational point on $C$. Substituting $s=3/10$ and $t=-7/10$ produces an elliptic curve whose LMFDB label is \texttt{150.c3}. Its mod-$30$ Galois image is conjugate to an index-$3$ subgroup of the reduction of $G^{\intercal}$ modulo $30$.

Write $G(30)$ for the reduction modulo $30$ of a suitable representative of $G$ in its conjugacy class. Reducing the adelic image generators listed for \href{https://www.lmfdb.org/EllipticCurve/Q/150a3/}{\texttt{150.c3}} in \cite{LMFDB} modulo $30$ gives
\[
H=\left\langle
\begin{psmallmatrix}1&20\\10&21\end{psmallmatrix},
\begin{psmallmatrix}11&16\\0&11\end{psmallmatrix},
\begin{psmallmatrix}26&5\\15&16\end{psmallmatrix},
\begin{psmallmatrix}6&5\\25&22\end{psmallmatrix}
\right\rangle.
\]
Direct computation gives $H\subseteq G(30)^{\intercal}$, $|H|=960$ and $|G(30)|=2880$, proving the asserted index.

\subsection{Acknowledgements} The authors of the appendix would like to thank LMFDB Modular Curves workshop held during March 2024 where the discussion on computing and simplifying universal elliptic curve equations started. The first author of appendix would also like to thank Jeremy Rouse for several valuable initial discussions.

\end{document}